\documentclass[12pt]{article}
\usepackage{enumerate}
\usepackage[square, comma, sort&compress, numbers]{natbib}
\usepackage{url} 

\newcommand{\blind}{1}

\RequirePackage{amsthm,amsmath,amsfonts,amssymb}

\usepackage{url}
\usepackage{bm}
\usepackage{graphicx}
\usepackage{color,xcolor}
\usepackage[ruled]{algorithm2e}
\usepackage{longtable}
\usepackage{supertabular}
\usepackage{booktabs}
\usepackage{multirow}
\usepackage{listings} 
\usepackage{xcolor} 
\usepackage{booktabs,dcolumn,caption}
\usepackage{lineno}
\usepackage{url}

\usepackage[ruled]{algorithm2e}
\def\be{\begin{equation}}
\def\ee{\end{equation}}
\def\bea{\begin{eqnarray}}
\def\eea{\end{eqnarray}}
\def\nn{\nonumber}
\def\fn{\footnote}

\newcommand{\vertiii}[1]{{\left\vert\kern-0.25ex\left\vert\kern-0.25ex\left\vert #1 
    \right\vert\kern-0.25ex\right\vert\kern-0.25ex\right\vert}}

\newtheorem{theorem}{Theorem}
\newtheorem{lemma}{Lemma}
\newtheorem{proposition}{Proposition}
\newtheorem{assumption}{Assumption}
\newtheorem{corollary}{Corollary}
\theoremstyle{remark}

\usepackage{setspace}
\doublespace 

\begin{document}          

\if1\blind
{
  \title{\bf Weighted Distributed Estimation under Heterogeneity}
  \author{Jia Gu\hspace{.2cm}\\
    Center for Statistical Science, Peking University, and  \\
    Song Xi Chen \thanks{corresponding author} \\
    School of Mathematical Science and Guanghua School of Management
\\ 
Peking University}
  \maketitle
} \fi

\if0\blind
{
  \bigskip
  \bigskip
  \bigskip
  \begin{center}
    {\LARGE\bf Weighted Distributed Estimation under Heterogeneity}
\end{center}
  \medskip
} \fi

\bigskip
\begin{abstract}
This paper considers distributed M-estimation under heterogeneous  distributions among distributed  data blocks. A weighted distributed  estimator is proposed to improve the efficiency of the standard "Split-And-Conquer" (SaC) estimator for the common parameter shared by all the data blocks.  
The weighted distributed estimator is shown to be at least as efficient as the would-be full sample  and the generalized method of moment estimators  with the latter two estimators  requiring full data access. 
 A bias reduction is formulated  to the WD estimator to accommodate  much larger numbers of data blocks than the existing methods without sacrificing the estimation efficiency, and a similar debiased operation is made to the SaC estimator. The mean squared error (MSE) bounds and the asymptotic distributions of the WD and the two debiased estimators are derived, which shows advantageous performance of the debiased estimators when the number of data blocks is large. 
\end{abstract}

\noindent%

\bigskip  

\section{Introduction} 
Modern big data have brought new challenges to statistical  inference. 
One such challenge is that despite the shear volume of the data, a full communication among the data points may not be possible 
 due to the expensive cost of data communication or the privacy concern. The 
 distributed  or the "split-and-conquer" (SaC) method  has been proposed to  divide the full data sample to smaller size data blocks as data communication is too costly to attain an estimation task in a timely fashion.  The SaC estimator is also suited to the situations where the data are naturally divided to data blocks and data communication among the data blocks are prohibited due to privacy concern.  
 The 
  SaC  estimation had been considered in \cite{Lin-Xi-2010} for the 
  U-statistics, \cite{Zhang-2013} for the M-estimation, \cite{Chen-Xie-2014} for the generalized linear models,  \cite{Volgushev-2017} and \cite{Chen-Liu-2019} for the quantile regression, while \cite{Battey-2018} studied the high dimensional testing and estimation with sparse penalties. 
\cite{Chen-Peng-2021} studied  the estimation efficiency and asymptotic  distributions for
  the general asymptotic symmetric statistics \citep{Lai-Wang-1993}  
and found differences in the efficiency and the asymptotic distributions between the non-degenerate and degenerate cases. 

Bootstrap resampling-based methods had  been introduced 
to 
facilitate  statistical inference. 
 \cite{Kleiner-2014} proposed the bag-of-little bootstrap (BLB) method for 
the  plug-in estimators by making up  economically 
the full sample for the  distributed inference. \cite{Sengupta-2015} suggested a sub-sampled double bootstrap method designed to improve the computational efficiency of the BLB. 
\cite{Chen-Peng-2021}  proposed the distributed and the pseudo-distributed bootstrap methods with the former conducted the resampling within each data block while the latter directly resampled the distributed  statistics. 

Privacy has been a major concern in big data applications where people are naturally reluctant to share the raw data to form a pool of big data as practised in  the traditional full sample estimation. However, the data holders may like to contribute summary statistics without having to give away the full data information.   Federated Learning 
 or the distributed inference with a central host has been proposed to accommodate such reality \cite{pmlr-v54-mcmahan17a, Yang-Liu-2019, Li-Smith-2020,Kairouz-2021},  
where  summary statistics of the data blocks or the gradients of the objective functions associated with the private data blocks 
 are submitted to a central host for forming aggregated estimation or computation.

Homogeneous distribution among the data blocks are assumed in majority of the distributed inference studies with only a few exceptions  \citep{Zhao-Cheng-2014,Duan-2021}. Federated Learning, on the other hand, was introduced to mitigate many challenges arising from classical distributed optimization. In particular, heterogeneous or Non-IID distributed data across different data blocks is one of the defining characteristics and challenges in the Federated Learning \cite{Li-Smith-2020,Kairouz-2021}. Indeed, it is natural to expect 
 the existence of heterogeneity,  
especially for data stored in different locations or generated by different stochastic mechanism, for instance mobile phones of different users. { However, there has been little published works on the statistical properties of estimators considered in the Federated Learning.} 

This paper considers distributed 
 estimation under heterogeneous  distributions among the data blocks, which is {closely related to}   the Federated Learning  and  especially  the multi-task learning (MTL) \citep{Zhang-Yang-IEEE-2021}. 
{ We consider distributed M-estimation} where there is a common parameter shared by the  distributions of the data blocks and  data-block specific heterogeneous parameters.  Our treatment of the heterogeneity is made by explicit parameterization, which is different from the MTL where the heterogeneity is regularized by penalty terms. 
%
It is noted that \cite{Duan-2021}  considered a heterogeneous setting, but under a fully parametric likelihood framework.
%
{ Our study reveals} that in the presence of the heterogeneity 
the full sample M-estimator of the common parameter obtained by requiring  full data communication,  
{ can be less} efficient  than the SaC estimator. However,  
this phenomenon disappears if the objective function  of the M-estimation  satisfies a generalized second-order Bartlett's identity, 
which are satisfied by the parametric and quasi likelihoods, and the least square estimation in the parametric regression.  

We propose a weighted distributed (WD) estimator, which is 
 {asymptotically at least  as efficient as} the full sample and the SaC estimator when the number of data blocks $K = o(N^{1/2})$ { where $N$ is the full sample size}. The mean-squared error (MSE) bound and the asymptotic distribution of the proposed WD estimator are derived, { as well as the asymptotic equivalence between the WD and the generalized method of moment (GMM) estimator. 
 We propose a debiased  weighted 
  distributed (dWD) estimator with  a data splitting mechanism on each data block to remove the correlation between the empirical bias correction and the weights used to tackle the heterogeneity.  {The dWD is asymptotically as efficient as the WD estimator, but with a more relaxed constraint of $K = o(N^{2/3})$.} 
  The bias-correction is also applied to the SaC formulation leading to a more communication-efficient 
  dSaC estimator, which is shown to be 
 more  accurate  than the subsampled average mixture estimator 
 (SAVGM) \cite{Zhang-2013} in the  homogeneous case. }   
 

The paper is organized as follows. The estimation framework and necessary notations for the study are outlined in Section \ref{sec: Preliminaries}. The relative efficiency between the full sample and the SaC estimators  under the heterogeneity is discussed in Section \ref{sec: full vs sac} to motivate the construction of the weighted distributed (WD) estimator. 
The WD estimator is introduced in Section \ref{sec: WD} along with its efficiency, asymptotic distribution and MSE bound. 
Statistical properties of two debiased estimators dSaC and dWD are revealed  in Section \ref{sec: bias reduction}. Section \ref{sec: Simulation} provides numerical verification to the theoretical results. Section \ref{sec: Discussion} concludes with a discussion. Technical details are reported in the supplementary materials (SM). 


\section{Preliminaries} \label{sec: Preliminaries}

Suppose that there is  a large data sample of size $N$, which is divided into $K$ data blocks of sizes  $\{n_k\}_{k=1}^K$ such that  
$N = \sum_{k=1}^Kn_k$ and  let $n = NK^{-1}$ be the average sample size of the data blocks. 
For the relative sample size among data blocks, we assume the following assumption.  
\begin{assumption}
\label{assumption: relative sample size}
There exist constants $ 0 < c < 1 < C$ such that $ c \leq\frac{n_{k_1}}{n_{k_2}}\leq C$ for all pairs of $(k_1,k_2)$, 
and if $K$ is a fixed constant we further assume that $\frac{n_k}{N}\rightarrow \gamma_k \in (0,1)$ for a set of constants $\{\gamma_k\}_{k=1}^K
$.
\end{assumption}


The $k$-th data block consists of a sub-sample  $ \{X_{k,i} \}_{i=1}^{n_k}$ which  are independent and identically distributed (IID) random vectors from a 
probability space $(\Omega,\mathcal{F},P)$ to $(\mathbb{R}^d,\mathcal{R}^d)$ with $F_{k}$ as the distribution.  
{The K distributions $\{F_k\}$ share a common parameter }
 $\phi \in \mathbb{R}^{p_1}$, while each $F_k$ 
 has another 
 parameter $\lambda_k \in \mathbb{R}^{p_2}$ 
 specific to $F_{k}$ of the $k$-th data block. There are maybe other hidden parameters which define $F_k$, which are however not directly involved in the semi-parametric  
 M-estimation, and thus are not of interest in the study.

 The parameters of interests in the $k$-th block are $\theta_k = (\phi^T,\lambda_k^T)^T$, 
 and the overall 
 parameters of interests  are $\theta = (\phi^T,\lambda_1^T,\lambda_2^T,...,\lambda_K^T)^T \in \mathbb{R}^{p_1 + Kp_2}$.
  Suppose there is a {common} 
  objective function $M(X;\phi,\lambda_k)$ { that is convex with respect to the parameter $(\phi,\lambda_k)$} and facilitates the M-estimation of the parameters {in each data block}. In general, the criteria function can be made block specific, say $M_k$ function.  Indeed,  the presence of the heterogeneous local parameters $\{\lambda_k\}_{k=1}^K$ leads to  different $M_k(x, \phi) = M(x, \phi, \lambda_k)$ for the inference on $\phi$,   which connects to the  multi-task learning (MTL). 
  
In the $k$-th data block the true parameter $\theta_k^*  = (\phi^{*T},\lambda_k^{*T})^T$ 
is defined as the unique minimum of the expected objective function, namely  
\begin{equation} 
\theta_k^*  = (\phi^{*T},\lambda_k^{*T})^T = \underset{\theta_k\in\Theta_k}{argmin}  \quad \mathbb{E}_{F_k} M(X_{k,1}; \phi, \lambda_k).  
\label{eq:local obj}
\end{equation}
 The  true common parameter $\phi^{\ast}$ appears in all $\theta_k^*$, and the block-specific  
 $\{\lambda_k^{*T}\}_{k=1}^K$ may differ from each other. 
 The entire set of true parameters 
 $\theta^* = (\phi^{*T},\lambda_1^{*T}, \cdots,\lambda_K^{*T})^T$,   
{
can  be also identified as 
\begin{equation}
    \theta^* = \underset{\theta \in \Theta}{argmin}\sum_{k=1}^K\gamma_k\mathbb{E}_{\theta_k^*}M(X_{k,1};\phi,\lambda_k). 
    \label{eq:global obj}
\end{equation} 

If the data could be shared across the data blocks,  we would attain the conventional full sample M-estimator  
\be 
\hat{\theta}_{full} = \underset{\theta \in \Theta}{argmin}\sum_{k=1}^K\sum_{i=1}^{n_k}M(X_{k,i};\phi,\lambda_k), \label{eq:full-Mestimation} 
\ee 
which serves as a benchmark for the distributed estimators. 
 Let $\psi_{\phi}(X_{k,i};\phi,\lambda_k) = \frac{\partial M(X_{k,i}; \phi,\lambda_k)}{\partial \phi}$ and $\psi_{\lambda}(X_{k,i};\phi,\lambda_k) = \frac{\partial M(X_{k,i};\phi,\lambda_k)}{\partial \lambda_k}$ 
be the score functions. { The estimating equations for the full sample M-estimators are } 
\begin{equation}
\begin{cases}
\sum_{k=1}^K\sum_{i=1}^{n_k}\psi_{\phi}(X_{k,i};\phi,\lambda_k) = 0, \\
\sum_{i=1}^{n_k}\psi_{\lambda}(X_{k,i};\phi,\lambda_k) = 0 \quad k = 1,...,K.
\end{cases}
\label{eq:full sample estimating equation}
\end{equation}  

{The above full sample estimation is not attainable for the distributed situations due to privacy or the costs associated with the data communications.} 
The distributed estimation first conducts local estimation on each data block, 
namely the local M-estimator
$$\hat{\theta}_k = (\hat{\phi}_k,\hat{\lambda}_k) = \underset{\theta_k\in\Theta_k}{argmin}  \sum_{i=1}^{n_k}M(X_{k,i};\theta_k)$$ 
with the corresponding estimating equations 
\begin{equation}
\begin{cases}
\sum_{i=1}^{n_k}\psi_{\phi}(X_{k,i};\phi_k,\lambda_k) = 0,\\
\sum_{i=1}^{n_k}\psi_{\lambda}(X_{k,i};\phi_k,\lambda_k) = 0. 
\end{cases}
\label{eq:local estimating equation}
\end{equation} 

Then, the "split-and-conquer" (SaC)  estimator for the common parameter $\phi$ is 
\begin{equation}
    \hat{\phi}^{SaC} = \frac{1}{N}\sum_{k=1}^Kn_k\hat{\phi}_k.
    \label{eq:SaC estimator}
\end{equation}

{The heterogeneity among the distributions  and the inference models 
among the data blocks bring new dimensions to the discussion of the relative efficiency and the estimation errors, which are the focus of this paper.  We are to show that the conventionally weighted SaC estimator 
(\ref{eq:SaC estimator}) may not be the best formulation for the estimation of 
 $\phi$. } 
{Throughout this paper,  unless otherwise stated, $\|\cdot\|_2$ and $\vertiii{\cdot}_2$ represent the $L_2$ norm of a vector and a matrix, respectively. Besides,                   we will use $C$ and $C_i$ to denote absolute positive constants independent of $(n_k,K,N)$.}


{An important question is the  efficiency and the estimation errors of the SaC estimator $\hat{\phi}^{SaC}$ relative to the 
full sample estimator $\hat{\phi}_{full}$. For the homogeneous case,}  Chen and Peng (2021) \cite{Chen-Peng-2021}  
found that for the asymptotic symmetric statistics, the SaC  estimator (\ref{eq:SaC estimator}) attains the same efficiency 
of the full sample estimator in the non-degenerate case, but encounters an efficiency loss in the degenerate case due to a lack of communications among different data blocks. 
Zhang et al. (2013) \cite{Zhang-2013} 
{ derived the mean square error (MSE) bound for the SaC estimator { in the homogeneous case} and  showed that whenever $K\leq \sqrt{N}$, the SaC estimator achieves the best possible rate of convergence when all $ N$ samples are accessible. }

Consider the simultaneous estimating equations  of the  full sample M-estimation
\be 
\Psi_N(\mathbf{X};\theta) = \begin{pmatrix}\sum_{k=1}^K\sum_{i=1}^{n_k}\psi_{\phi}(X_{k,i};\phi,\lambda_k) \\ \sum_{i=1}^{n_1}\psi_{\lambda}(X_{1,i};\phi,\lambda_1)\\ \vdots \\ 
\sum_{i=1}^{n_K}\psi_{\lambda}(X_{K,i};\phi,\lambda_K)\end{pmatrix}. \label{eq: FOC of full objective} 
\ee
Define 
\bea
\Psi_{\theta}(\theta_k) &=& (\Psi_{\phi}(\theta_k)^T,\Psi_{\lambda}(\theta_k)^T)^T = \mathbb{E}\nabla_{\theta_k}M(X_{k,1};\theta_k),\nn
\\
\Psi_{\theta}^{\theta}(\theta_k) &=& \begin{pmatrix} \Psi_{\phi}^{\phi}(\theta_k) & \Psi_{\phi}^{\lambda}(\theta_k)\\
\Psi_{\lambda}^{\phi}(\theta_k) & \Psi_{\lambda}^{\lambda}(\theta_k)\end{pmatrix} = \mathbb{E}\nabla_{\theta_k}^2M(X_{k,1};\theta_k)
, \nn\\
J_{\phi|\lambda}(\theta_k) &=& \Psi_{\phi}^{\phi}(\theta_k) - \Psi_{\phi}^{\lambda}(\theta_k)\Psi_{\lambda}^{\lambda}(\theta_k)^{-1}\Psi_{\lambda}^{\phi}(\theta_k) \quad \hbox{and} \nn\\
S_{\phi}(X_{k,i};\theta_k) &=& \psi_{\phi}(X_{k,i};\theta_k) - \Psi_{\phi}^{\lambda}(\theta_k)\Psi_{\lambda}^{\lambda}(\theta_k)^{-1}\psi_{\lambda}(X_{k,i};\theta_k).\nn
\eea
Then we can apply Taylor's expansion and obtain (see Section 1.1 in SM for details)
\be
         \hat{\phi}_{full} - \phi^* =  -\{\sum_{k=1}^K\frac{n_k}{N}J_{\phi|\lambda}(\theta_k^*)\}^{-1}\frac{1}{N}\big\{\sum_{k=1}^{K}\sum_{i=1}^{n_k}S_{\phi}(X_{k,i};\theta_k^*)\big\} + o_p(N^{-1/2}),
    \label{eq:full expansion}
\ee


For the   local  estimator $(\hat{\phi}_k,\hat{\lambda}_k)$ based on the $k$-th  data block that solves  (\ref{eq:local estimating equation}), by replicating the same derivation leading to  (\ref{eq:full expansion}), we have 
\be
    \begin{cases}
         \hat{\phi}_k - \phi^* &= -n_k^{-1}J_{\phi|\lambda}(\theta_k^*)^{-1}\sum_{i=1}^{n_k}S_{\phi}(X_{k,i};\theta_k^*) + o_p(n_k^{-1/2}),\\
         \hat{\lambda}_k - \lambda_k^* &= -n_k^{-1}J_{\lambda|\phi}(\theta_k^*)^{-1}\sum_{i=1}^{n_k}S_{\lambda}(X_{k,i};\theta_k^*) + o_p(n_k^{-1/2}),
    \end{cases}
\ee
where 
\bea
J_{\lambda|\phi}(\theta_k) &=& \Psi_{\lambda}^{\lambda}(\theta_k) - \Psi_{\lambda}^{\phi}(\theta_k)\Psi_{\phi}^{\phi}(\theta_k)^{-1}\Psi_{\phi}^{\lambda}(\theta_k) \quad \hbox{ and }\nn\\
S_{\lambda}(X_{k,i};\theta_k) &=& \psi_{\lambda}(X_{k,i};\theta_k) - \Psi_{\lambda}^{\phi}(\theta_k)\Psi_{\phi}^{\phi}(\theta_k)^{-1}\psi_{\phi}(X_{k,i};\theta_k).\label{eq: score for lambda}
\eea

{The distributed inference setting is closely related to  the Multi-Task Learning (MTL) which fits separate local parameters 
 $\phi_k \in \mathbb{R}^p$ to the data of different data blocks (tasks) through 
 convex loss functions $\{\ell_k\}$. {In particular, the} MTL is formulated as 
 \cite{Smith-NIPS-2017}: 
  \be
        \underset{\Phi, \Omega}{min}\quad \bigg\{\sum_{k=1}^K\sum_{i=1}^{n_k}\ell_k(\phi_k^TX_{k,i},Y_{k,i}) + \mathcal{R}(\Phi, 
        \Omega)\bigg\},\label{eq: MTL}
        \ee
where $\{(X_{k,i},Y_{k,i}), i = 1,2,\cdots,n_k\}$
are data in the $k$-th block, {$\Phi$ is the matrix with $\{\phi_k\}_{k=1}^K$ as column vectors}, $\Omega \in \mathbb{R}^{K\times K}$ 
and $\mathcal{R}(\cdot,\cdot)$ measures  the extent of the heterogeneity among different data blocks. {Choices of  
$\mathcal{R}(\cdot,\cdot)$ include} 
$\mathcal{R}(\Phi,\Omega) = \delta_1 tr(\Phi\Omega \Phi^T) + \delta_2 \|\Phi\|_F^2$
for $\delta_1,\delta_2 >0$ and $\Omega = I_{K\times K} - \frac{1}{K}1_K1_K^T$ such that  $tr(\Phi\Omega\Phi^T) = \sum_{k=1}^K\|\phi_k - \bar{\phi}_K\|_2^2$ where $   
 \bar{\phi}_K = \frac{1}{K}\sum_{k=1}^K\phi_k$, which leads to 
 the  mean-regularized MTL \cite{Evgeniou-2004}.  
The second term of $\mathcal{R}$ performs 
regularization on each local model, trying to control the magnitude of the estimates of $\phi_k$. 

{The  distributed framework is well connected to the MTL in two key aspects. One is that despite we use the same objective (loss) function  $M$ over the data blocks, the  heterogeneity induced by local parameters $\{\lambda_k\}_{k=1}^K$ and the distributions  effectively define $M_k(\phi, x) = M(x, \phi, \lambda_k)$, which is equivalent to the block specific loss functions $\ell_k$ used in MTL. Another aspect is that although the MTL assumes different parameters $\{\phi_k\}$ over the data blocks, it regularizes them toward a common one. In contrast, we assume there is a common parameter $\phi$ shared by the heterogeneous distributions. }

\section{Full Sample versus SaC Estimation}\label{sec: full vs sac}

It is naturally expected that the full sample estimator $\hat{\phi}_{full}$ should be { at least as efficient as} the distributed  SaC estimator $\hat{\phi}^{SaC}$ since the former 
utilizes the full sample { information including the communications among different data blocks}. However, { we are to show that} this is not necessarily true in the presence of heterogeneity. To appreciate this point, 
we first list more 
 regularity conditions needed in the analysis. 

\begin{assumption}\label{assumption:wellsep}(\textbf{Identifiability})
The parameters $\theta_k^* = (\phi^*,\lambda_k^*) $ is the unique minimizer of $M_k(\theta_k)= \mathbb{E}M(X_{k,1};\theta_k)$ for $\theta_k \in \Theta_k$.
\end{assumption}
\begin{assumption}\label{assumption:compactness}(\textbf{Compactness}) The parameter space $\Theta_k$ is a compact and convex set in $\mathbb{R}^p$ and the true parameter $\theta_k^*$ is an interior point of $\Theta_k$ and $\underset{\theta_k \in \Theta_k}{sup}\|\theta_k - \theta_k^*\|_2\leq r$ for all $k\geq 1$ and some $r > 0$. The true common parameter $\phi^*$ is an interior point of a compact and convex set $\Phi\subset \Theta_k$.
\end{assumption}
\begin{assumption}
\label{assumption: Local strong convexity}(\textbf{Local strong convexity})
The population objective  function on the k-th data block $M_k(\theta_k) = \mathbb{E}M(X_{k,1};\theta_k)$ is twice differentiable, and there exists a constant $\rho_{-} > 0$ such that $\nabla^2_{\theta_k}M_k(\theta_k^*)\succeq \rho_{-}I_{p\times p}$. Here $A  \succeq B$ means $A - B$ is a positive semi-definite matrix.
\end{assumption}

{ These three assumptions are standard ones on the parameter space 
 and population objective functions as those in Zhang et al. (2013) \cite{Zhang-2013} and Jordan et al. (2019) \cite{Jordan-2019} for the homogeneous case. In the heterogeneous setting, Duan et al. (2021) \cite{Duan-2021} only requires the parameter space for the common parameter to be bounded, i.e. $\|\phi - \phi^*\|\leq r$ under a fully parametric setting, while in our assumption, we need the overall parameter space to be bounded. This stronger assumption is needed since we do not fully specify the distributions $\{F_k\}_{k=1}^K$ of the random variables and will be useful when we derive the MSE bound for the 
  {weighted distributed}  estimator which will be proposed in Section \ref{sec: WD}. } 

\begin{assumption}(\textbf{Smoothness})\label{assumption: smoothness}
There are finite positive constants $R, L, v$ and $ v_1$ such that for all $k\geq 1$,
$\mathbb{E}\|\nabla_{\theta_k}M(X_{k,1};\theta_k^*)\|_2^{2v_1}\leq R^{2v_1} \text{ and } \mathbb{E}\vertiii{\nabla_{\theta_k}^2M(X_{k,1};\theta_k^*) - \nabla_{\theta_k}^2M_k(\theta_k^*)}_2^{2v}\leq L^{2v}.$  
In addition, for any $x\in\mathbb{R}^d$, $\nabla_{\theta_k}^2M(x;\theta_k)$ and $\nabla_{\theta_k}M(x;\theta_k)\nabla_{\theta_k}M(x;\theta_k)^T$ are $G(x)-$ and $B(x)-$Lipschitz continuous, respectively, in the sense that
\begin{small}
$$\vertiii{\nabla_{\theta_k}^2M(x;\theta_k) - \nabla_{\theta_k}^2M(x;\theta_k^{'})}_2\leq G(x)\|\theta_k - \theta_k^{'}\|_2,$$ $$\vertiii{\nabla_{\theta_k}M(x;\theta_k)\nabla_{\theta_k}M(x;\theta_k)^T - \nabla_{\theta_k}M(x;\theta_k^{'})\nabla_{\theta_k}M(x;\theta_k^{'})^T}_2\leq B(x)\|\theta_k - \theta_k^{'}\|_2,$$
\end{small}
for all $\theta_k,\theta_k^{'}\in U_k := \{\theta_k|\|\theta_k - \theta_k^*\|_2 \leq \rho\}$ for some $\rho >0$, and $\mathbb{E}G(X_{k,1})^{2v} \leq G^{2v}, \mathbb{E}B(X_{k,1})^{2v}\leq B^{2v}$ for some positive constants $G$ and $B$. 
\end{assumption}

{The Lipschitz continuity of the outer product of the first-order derivative is required to control the estimation error when we estimate the asymptotic covariance matrix of the local estimator $\hat{\theta}_k$, and it can be directly verified under the logistic regression case; see Section 1.2 in the SM for details.}  

\begin{proposition}
Under Assumptions \ref{assumption: relative sample size} -  \ref{assumption: Local strong convexity} and Assumption \ref{assumption: smoothness} with $v, v_1 \geq  1$,
and if $K$ is { fixed}, 
then 
$\hat{\theta}_k\overset{P}{\rightarrow}\theta_k^*$ and $\hat{\theta}_{full} \overset{P}{\rightarrow}\theta^*$; 
$\hat{\phi}^{SaC} = \frac{1}{N}\sum_{k=1}^Kn_k\hat{\phi}_k$ and $\hat{\phi}_{full}$ are consistent to $\phi^{\ast}$. 
\end{proposition}

\begin{theorem}\label{thm: Compare SaC and full}
Under Assumptions \ref{assumption: relative sample size} -  \ref{assumption: Local strong convexity} and Assumption \ref{assumption: smoothness} with $v, v_1\geq  2$,  if $K$ is a fixed constant, the SaC estimator $\hat{\phi}^{SaC}$ and the full sample estimator $\hat{\phi}_{full}$ satisfy 
\begin{small}
\begin{subequations}
\begin{align}
    \sqrt{N}(\hat{\phi}^{SaC} - \phi^*) &\overset{d}{\rightarrow}\mathcal{N}(\mathbf{0},\sum_{k=1}^K\gamma_kJ_{\phi|\lambda}(\theta_k^*)^{-1}\Sigma_k(\theta_k^*)J_{\phi|\lambda}(\theta_k^*)^{-1}),\\
\sqrt{N}(\hat{\phi}_{full} - \phi^*) &\overset{d}{\rightarrow}\mathcal{N}(\mathbf{0},(\sum_{k=1}^K\gamma_kJ_{\phi|\lambda}(\theta_k^*))^{-1}(\sum_{k=1}^K\gamma_k\Sigma_k(\theta_k^*))(\sum_{k=1}^K\gamma_kJ_{\phi|\lambda}(\theta_k^*))^{-1}),
\end{align}
\end{subequations}
where $J_{\lambda|\phi}(\theta_k^*) = \Psi_{\lambda}^{\lambda}(\theta_k^*) - \Psi_{\lambda}^{\phi}(\theta_k^*)\Psi_{\phi}^{\phi}(\theta_k^*)^{-1}\Psi_{\phi}^{\lambda}(\theta_k^*)$ and $ \Sigma_k = Var\{S_{\phi}(X_{k,1};\theta_k^*)\}$.
\end{small}
\end{theorem} 

Define $V(\Sigma,A) = (A^{T})^{-1}\Sigma A^{-1}$ as a mapping from $\mathbb{S}_{++}^{p_1\times p_1} \times GL(\mathbb{R}^{p_1})$ to $\mathbb{S}_{++}^{p_1\times p_1}$, where $\mathbb{S}_{++}^{p_1\times p_1} $ and $GL(\mathbb{R}^{p_1})$ denote the symmetric positive definite matrices and invertible real matrices of order $p_1$, respectively. Since $\Sigma_{k=1}^K\gamma_k = 1$ and $\gamma_k > 0$, the asymptotic variance of $\hat{\phi}^{SaC}$ can be interpreted as a convex combination of function values $\{V(\Sigma_k(\theta_k^*),J_{\phi|\lambda}(\theta_k^*))\}_{k=1}^K$ and that of $\hat{\phi}_{full}$  can be seen as $V(\sum_{k=1}^K\gamma_k\Sigma_k(\theta_k^*),\sum_{k=1}^K\gamma_kJ_{\phi|\lambda}(\theta_k^*))$. However, $V(\cdot,\cdot)$ is not  convex with respect to its arguments $(\Sigma, A)$, which means that the inequality 
\be
\{\sum_{k=1}^K\gamma_kJ_{\phi|\lambda}(\theta_k^*)\}^{-1}\{\sum_{k=1}^K\gamma_k\Sigma_k(\theta_k^*)\}\{\sum_{k=1}^K\gamma_kJ_{\phi|\lambda}(\theta_k^*)\}^{-1} \preceq \sum_{k=1}^K\gamma_kJ_{\phi|\lambda}(\theta_k^*)^{-1}\Sigma_k(\theta_k^*)J_{\phi|\lambda}(\theta_k^*)^{-1}\nn
\ee
does not always hold. In other words, $\hat{\phi}_{full}$ is not necessarily more efficient than $\hat{\phi}^{SaC}$.

To gain 
understanding of Theorem \ref{thm: Compare SaC and full} and to motivate the  
 {weighted distributed} estimator, 
  we consider 
the {errors-in-variables model}. 
Suppose that one observes $K$ blocks of independent data samples $\{(X_{k,i},Y_{k,i})\}_{i=1}^{n}$ for $k = 1,2...,K$ and $N = nK$, where $(X_{k,i},Y_{k,i})$ are IID and generated from 
 the following model: 
\begin{equation}
    \begin{cases}
X_k = Z_k + e_k,\\
Y_k = \phi^* + \lambda_k^* Z_k + f_k,\\
\end{cases}
\label{eq:error-in-variable} 
\end{equation}
where $\{Z_k\}_{k=1}^K$ are random variables whose measurements $\{(X_k,Y_k)\}_{k=1}^K$ are subject to errors $\{(e_k,f_k)\}_{k=1}^K$, and $(e,f)$ is bivariate normally distributed with zero mean and covariance matrix $\sigma^2I_2$ and is independent of $Z_k$. Here, 
$\phi^*$ is the common parameter across all  data blocks while $\lambda_k^*  (\lambda_k^* > 0)$ represents the block specific parameter.  We assume that $Var(e) = Var(f)$ to avoid any identification issue arisen when $Z$ is also normally distributed \cite{Reiersol-1950}. There is a considerable literature on the regression problem with measurement errors,  as  summarised in \cite{Fuller-1987,Schafer-1996}.  

We consider the approach displayed in Example 5.26 of \cite{Vaart-1999}  which  constructs a kind of marginal likelihood followed by  centering to make a bona fide score equation, as detailed in Section 1.3 of the SM. The M-function is 
\be
M(X_k,\theta_k) = \frac{1}{2\sigma^2(1 + \lambda_k^2)}(\lambda_kX_k - (Y_k - \phi))^2, \label{eq: obj for EIV model}
\ee
{ with the score equation satisfying } 
$\mathbb{E}\nabla M(X_{k,1},Y_{k,1}|Z_{k,1},\theta_k^*) = \mathbf{0}_{2\times 1}$. 

For simplicity we assume $K = 2$, then from Theorem \ref{thm: Compare SaC and full} we have
\be\label{eq: compare SaC and full in EIV}
\begin{cases}
Var(\hat{\phi}_{full}) \approx \{\frac{\sigma^2\mathbb{E}Z^2}{var(Z)} \frac{2}{\frac{1}{1+\lambda_1^{*2}}+\frac{1}{1+\lambda_2^{*2}}} + \frac{\sigma^4(\mathbb{E}Z)^2}{var^2(Z)}\frac{\frac{2}{(1+\lambda_1^{*2})^2}+\frac{2}{(1+\lambda_2^{*2})^2}}{(\frac{1}{1+\lambda_1^{*2}}+\frac{1}{1+\lambda_2^{*2}})^2}\}\frac{1}{N},\\

Var(\hat{\phi}^{SaC}) \approx \{\frac{\sigma^2\mathbb{E}Z^2}{var(Z)}\frac{(1+\lambda_1^{*2})+(1+\lambda_2^{*2})}{2} + \frac{\sigma^4(\mathbb{E}Z)^2}{var^2(Z)}\}\frac{1}{N}.
\end{cases}
\ee
Note that { the coefficients 
to $\frac{\sigma^2\mathbb{E}Z^2}{var(Z)}$ in the first terms of the variances } 
are harmonic  and arithmetic means of $\{1+\lambda_1^{*2},1+\lambda_2^{*2}\}$, respectively. By the mean inequality the coefficient 
in  the first term of $Var(\hat{\phi}^{SaC})$ is larger than  that 
in  $Var(\hat{\phi}_{full})$.  The second term of the variances involves   
$(\mathbb{E}Z)^2$  as a multiplicative factor. Thus, if the unobserved  $Z$ has zero mean, 
the full-sample estimator would be at least as good as the SaC estimator in terms of variation when the sample size goes to infinity.  However, the story may change 
when $\mathbb{E}Z \not = 0$, because the second term of  $Var(\hat{\phi}_{full})$ { has a factor} which is 
the square 
of a ratio between the quadratic mean and the arithmetic mean of $(\frac{1}{1+\lambda_1^{*2}},\frac{1}{1+\lambda_2^{*2}})$. The factor is larger than  { or equal to 1 if and only if  $\lambda_1^* = \lambda_2^*$ namely the homogeneous case. } 
In the heterogeneous case, by adjusting  $\frac{\sigma^4(\mathbb{E}Z)^2}{var^2(Z)}/\frac{\sigma^2\mathbb{E}Z^2}{var(Z)} $, we can  find cases such that $\lambda_1^* \ne \lambda_2^*$
such that 
the full sample estimator { has a larger variance than} the SaC estimator. Simulation experiments presented in Section \ref{sec: Simulation} display such cases. 

\section{Weighted Distributed Estimator} 
\label{sec: WD}
The previous section shows that the full sample estimator $\hat{\phi}_{full}$ under heterogeneity may be less  efficient than the simple averaged $\hat{\phi}^{SaC}$. 
This phenomenon suggests that { the conventional wisdom in the homogeneous 
context case may not be applicable to the heterogeneous case.  One may also wonder if the simple SaC estimator can be improved under the heterogeneity.  Specifically, how to better aggregate the local estimator $\hat{\phi}_k$  
for more efficiency estimation to the common parameter $\phi$  is the focus of  this section. } 

\subsection{Formulation and Results} 

Consider a class of estimators formed by linear combinations of the local estimators  $\{\hat{\phi}_k\}$: 
$$\{\hat{\phi}_{w}^{SaC}| \hat{\phi}_{w}^{SaC} = \sum_{k=1}^KW_k\hat{\phi}_k,W_k\in\mathbb{R}^{p_1\times p_1},\sum_{k=1}^KW_k = I_{p_1}\}.$$ We want to minimize the asymptotic variance of $\hat{\phi}_{w}^{SaC}$ with respect to $\{W_k\}_{k=1}^K$. According to a generalization of Theorem \ref{thm: Compare SaC and full} 
\be
    AsyVar(\hat{\phi}_{w}^{SaC}) = \sum_{k=1}^Kn_k^{-1}W_kA_k^{-1}\Sigma_k(A_k^T)^{-1}W_k^T,
\ee
 where $A_k = J_{\phi|\lambda}(\theta_k^*)$ and $\Sigma_k = Var\{S_{\phi}(X_{k,i};\theta_k^*)\}$. 
 {It is noted that the asymptotic variance is defined via the asymptotic normality of the M-estimation.} For the { time being},  $A_k$ and $\Sigma_k$ are assumed known and denote $H_k = A_k^{-1}\Sigma_k(A_k^T)^{-1}$.  We choose the trace operator 
{ as a measure on}  the size of the asymptotic covariance matrix and this leads to the minimization problem
\be
\underset{W_k}{Minimize} \quad tr\bigg(\sum_{k=1}^Kn_k^{-1}W_kH_kW_k^T\bigg) \quad 
s.t. \quad \sum_{k=1}^K W_k = I_{p_1}, 
\label{eq:Best weight}
\ee
which is a 
 convex optimization problem. It can be solved via the  Lagrangian multiplier method which gives $W_k^{*} = (\sum_{s=1}^K n_sH_s^{-1})^{-1}n_kH_k^{-1} $. {If we replace the trace with the Frobenius norm in  the objective function (\ref{eq:Best weight}), the same solution is attained as shown in Section 1.4 of the SM.} 
 The SaC estimator under the optimal weights $W_k^*$ is called the  weighted distributed (WD) estimator  and denoted as $\hat{\phi}^{WD}$. By construction, the WD estimator is at least as efficient as 
 the SaC estimator (\ref{eq:SaC estimator}).  To compare the relative efficiency between $\hat{\phi}_{full}$ and $\hat{\phi}^{WD}$, we note that
\bea 
AsyVar(\hat{\phi}_{full}) &=& \big\{(\sum_{k=1}^Kn_kA_k)^T(\sum_{k=1}^Kn_k\Sigma_k)^{-1}(\sum_{k=1}^Kn_kA_k)\big\}^{-1} \quad \hbox{and}  \nn \\
AsyVar(\hat{\phi}^{WD}) &=& \left(\sum_{k=1}^Kn_kA_k^T\Sigma_k^{-1}A_k\right)^{-1}.\label{eq: compare var between full and WD}
\eea 
 
Define $F(\Sigma,A) = A^T\Sigma^{-1}A$. 
If we can show the convexity of $F$,  an application of Jensen's inequality will establish the relative efficiency of the two estimators. In fact, we have the following lemma. 

\begin{lemma}\label{lemma: sandwich inequality}
Suppose $H$ and $K$ are positive definite matrices of order $p$, and $X$ and $Y$ are arbitrary $p\times m$ matrices.  Then, 
\be
Q = X^TH^{-1}X + Y^TK^{-1}Y - (X+Y)^T(H+K)^{-1}(X+Y)\succeq \mathbf{0}.\nn
\ee
\label{lemma:convexity}
\end{lemma}
The lemma implies that 
$$(\sum_{k=1}^Kn_kA_k)^T(\sum_{k=1}^Kn_k\Sigma_k)^{-1}(\sum_{k=1}^Kn_kA_k) \preceq \sum_{k=1}^Kn_kA_k^T\Sigma_k^{-1}A_k,$$
which means that the WD estimator is at least as efficient as the full sample estimator, and can be more efficient than $\hat{\phi}_{full}$. {That is to say, the simultaneous estimating equations (\ref{eq: FOC of full objective}), which are obtained from the first-order derivative of the the simple summation of local objectives $\sum_{i=1}^{n_k}M(X_{k,i};\theta_k)$, are not the best formulation of the M-estimation problem, since the formulation itself does not utilize the heterogeneity existed in the data blocks. In contrast, the WD estimator exploits the potential efficiency gain from the heterogeneity by re-weighting of the local estimators, and this is why the full sample estimator may not be as efficient as the WD estimator.}

\subsection{Likelihood and Quasi-likelihood}

{The above results lead us to wonder } 
whether we can attain more efficient distributed estimators than the full sample estimator under the heterogeneity if we restrict to a fully parametric setting. 
When the distribution of $X_{k,i}$ is fully  parametric 
with density function  $f(\cdot;\phi,\lambda_k)$,   the 
Fisher information matrix in the $k$-th data block is 
$$I(\theta_k) = I(\phi,\lambda_k) = \begin{pmatrix}I_{\phi\phi}  & I_{\phi\lambda_k} \\
I_{\lambda_k \phi} & I_{\lambda_k\lambda_k}\end{pmatrix} = -\mathbb{E}\begin{pmatrix} 
\frac{\partial^2}{\partial \phi^2}logf(X_{k,1};\theta_k) & \frac{\partial^2}{\partial \phi\partial \lambda^T}logf(X_{k,1};\theta_k)\\
\frac{\partial^2}{\partial \lambda \partial \phi^T}logf(X_{k,1};\theta_k) & \frac{\partial^2}{\partial \lambda^2}logf(X_{k,1};\theta_k)
\end{pmatrix},$$
and the partial information matrix as $I_{\phi|\lambda_k} = I_{\phi\phi} - I_{\phi\lambda_k}I_{\lambda_k\lambda_k}^{-1} I_{\lambda_k\phi}$. 
Now, the objective function for the M-estimation (also the maximum likelihood estimation (MLE)) is $M(X_{k,i};\phi,\lambda_k) = -\log f(X_{k,i};\phi,\lambda_k)$. 
Routine 
derivations show that $\Sigma_k = Var\{S_{\phi}(X_{k,1};\theta_k^*)\}= I_{\phi|\lambda_k}$ and $A_k = J_{\phi|\lambda}(\theta_k^*) =I_{\phi|\lambda_k}$. Thus, 
\bea 
AsyVar(\hat{\phi}_{full}) &=& AsyVar(\hat{\phi}^{WD})   = \left(\sum_{k=1}^K n_kI_{\phi|\lambda_k}\right)^{-1} 
\quad \hbox{ and } \nn \\ 
AsyVar(\hat{\phi}^{SaC}) &=& \frac{1}{N^2}\sum_{k=1}^K n_kI_{\phi|\lambda_k}^{-1}.  \nn 
\eea 

{ A direct application of Lemma \ref{lemma: sandwich inequality} shows that}
\be
AsyVar(\hat{\phi}_{full}) = AsyVar(\hat{\phi}^{WD})\preceq AsyVar(\hat{\phi}^{SaC}). \label{eq: compare full, WD, SaC}
\ee
Thus,  the full sample MLE can automatically adjust for the heterogeneity 
 and is at least as efficient as SaC estimator $\hat{\phi}^{SaC}$. Besides, the weighted  distributed estimators  $\hat{\phi}^{WD}$  can fully recover the efficiency gap of the SaC estimator.

The same relationship among $\hat{\phi}_{full},\hat{\phi}^{SaC}$ and $\hat{\phi}^{WD}$ also holds for the maximum quasi-likelihood estimator (MQLE) with independent observations (see Section 1.5 in the SM for details).
If one looks into the asymptotic variances of the MLE and MQLE,  it  can be found that 
 the underlying reason for (\ref{eq: compare full, WD, SaC}) is that the  two special M-estimation functions satisfy the second order Bartlett's identity \cite{Bartlett-1953,McCullagh-1983}: 
$$\mathbb{E}\nabla M(X_k,\theta_k^*)\nabla M(X_k,\theta_k^*)^T =   \mathbb{E}\nabla^2 M(X_k,\theta_k^*). $$
By the variance formula of the asymptotic distribution of the M-estimator and Lemma \ref{lemma: sandwich inequality}, we readily have that the Bartlett's identity can be relaxed by inserting  a 
factor $\gamma \not = 0$ such that 
\be
\mathbb{E}\nabla M(X_k,\theta_k^*)\nabla^TM(X_k,\theta_k^*) = \gamma \mathbb{E}\nabla^2M(X_k,\theta_k^*).  \label{eq: full sample estimator efficiency}
\ee
An important example for such a case is the least square estimation for {the parametric} regression with homoscedastic and non-autocorrelated disturbances ({ see Section 1.6 in the SM for details}).  Otherwise the full sample  least square estimator may  not be efficient and there is an opportunity for the weighted distributed least square estimation. 
In summary, as long as the objective  function $M(x_k,\theta_k)$  
satisfies (\ref{eq: full sample estimator efficiency}), then $\hat{\phi}_{full}$ attains the same asymptotic efficiency as $\hat{\phi}^{WD}$,  and $\hat{\phi}^{SaC}$ is at most as efficient as the former two estimators.

\subsection{Relative to Generalized Method of Moment Estimation}  
{To provide a benchmark on the efficiency of the WD 
 estimation, we consider the generalized method of moment (GMM) estimator  \citep{Hansen-1982}. The GMM estimator possess certain optimal property 
for semiparametric inference that the weighted distributed estimation can compare with, despite the GMM requires more data sharing than the distributed inference would require.  }  
%

The score functions of the M-estimation on each data block can be aggregated and combined to form the  moment equations
\be
\begin{cases}
\sum_{i=1}^{n_k}\psi_{\phi}(X_{k,i};\phi,\lambda_k) = 0, \\
\sum_{i=1}^{n_k}\psi_{\lambda}(X_{k,i};\phi,\lambda_k) = 0, \quad k = 1,...,K.\label{eq: GMM}
\end{cases}
\ee
 There are $pK$ estimating equations, where  the dimension of 
 $\theta^*$ is $pK - (K-1)p_1$.  
 Thus,  
the parameter is over-identified which offers potential in efficiency gain for the GMM \citep{Hansen-1982}. 
The GMM estimation based on the  moment restrictions  (\ref{eq: GMM}) 
is asymptotically equivalent to 
solving  the following 
problem: 
\be
\hat{\theta}_{GMM} = \underset{\theta_k = (\phi,\lambda_k)\in\Theta_k, 1\leq k \leq K}{argmin}\quad \tilde{\psi}_N^T(\theta)W_0\tilde{\psi}_N(\theta),\label{eq: GMM objective}
\ee
where $W_0 = Var(\tilde{\psi}_N(\theta^*))^{-1}$ is the optimal weighting matrix \citep{Hansen-1982,Yaron-Hansen-1996-JBES} and 
\be
\tilde{\psi}_N(\theta) = (\sum_{i=1}^{n_1}\psi_{\phi}(X_{1,i};\theta_1)^T,\sum_{i=1}^{n_1}\psi_{\lambda}(X_{1,i};\theta_1)^T,\cdots,\sum_{i=1}^{n_K}\psi_{\phi}(X_{K,i};\theta_K)^T,\sum_{i=1}^{n_K}\psi_{\lambda}(X_{K,i};\theta_K)^T)^T\nn 
\ee

The asymptotic variance of the GMM estimator 
\citep{Hansen-1982} is  $AsyVar(\hat{\theta}_{GMM}) = (G_0^TW_0G_0)^{-1}$, where $G_0^T = \mathbb{E}\{\frac{\partial \tilde{\psi}_N^T(\theta^*)}{\partial \theta}\}$. 
A  derivation given in  Section 1.7 of the SM shows that 
\be
AsyVar(\hat{\phi}_{GMM}) =\{\sum_{k=1}^K n_kJ_{\phi|\lambda}\Sigma_k^{-1}J_{\phi|\lambda}\}^{-1}.
\ee
Thus, { the weighted distributed estimator's  efficiency is the same as that of} the GMM estimator. 
{ This is very encouraging to the proposed WD estimator  as  it attains the same efficiency as the GMM without requiring much data sharing among the  blocks, which avoids  the expenses of the data transmission and preserves the privacy of the data.} 

\subsection{Estimation of Weights in one round communication} 
\label{subsection: weights}

To formulate the WD estimator,  
 the optimal weights $W_k^* = (\sum_{s=1}^Kn_sH_s^{-1})^{-1}n_kH_k^{-1}$ have to be estimated. By the structure of $W_k^*$, we only need to separately estimate $H_k$, the leading principal submatrix of order $p_1$ of the asymptotic covariance matrix $\tilde{H}_k$ of $\hat{\theta}_k$. It is noted that 
\be
\tilde{H}_k = (\nabla \Psi_{\theta}(\theta_k^*))^{-1}\mathbb{E}\{\psi_{\theta_k}(X_{k,1};\theta_k^*)\psi_{\theta_k}(X_{k,1};\theta_k^*)^T \}(\nabla \Psi_{\theta}(\theta_k^*))^{-1}  = \begin{pmatrix} H_k & *\\ * & *\\ \end{pmatrix}\nn,
\ee
where $\Psi_{\theta}(\theta_k) = \mathbb{E}\psi_{\theta_k}(X_{k,1};\theta_k)$. We can construct  the  sandwich estimator \cite{Stefanski-Boos-2002-AmericanStatistician} to estimate $\tilde{H}_k$ and then $H_k$.   The {distributive} procedure to attain the WD estimator is  summarized in the Algorithm 1.
\begin{algorithm}
\caption{\textbf{W}eighted \textbf{D}istributed M-estimator }
\LinesNumbered
\KwIn{$\{X_{k,i},k = 1,...,K; i = 1,...,n_k\}$}
\KwOut{$\hat{\phi}^{WD},\hat{\lambda}_k$}
Obtain the initial estimates $\hat{\theta}_k = (\hat{\phi}_k,\hat{\lambda}_k)$ based on data block $k$ \;
Calculate $\hat{H}_k(\hat{\theta}_k)$ in each block, which is the leading principal  sub-matrix of order $p_1$ of 
$(\nabla_{\theta_k}\hat{\Psi}_{\theta_k})^{-1}(n_k^{-1}\sum_{i=1}^{n_k}\psi_{\theta_k}(X_{k,i};\hat{\theta}_{k})\psi_{\theta_k}(X_{k,i};\hat{\theta}_{k})^T)(\nabla_{\theta_k}\hat{\Psi}_{\theta_k})^{-T}$ where $\hat{\Psi}_{\theta_k} = n_k^{-1}\sum_{i=1}^{n_k}\psi_{\theta_k}(X_{k,i};\hat{\theta}_{k})$\;
Send $(\hat{\phi}_k,\hat{H}_k(\hat{\theta}_k)^{-1})$ to a central server and construct $\tilde{\hat{\phi}}^{WD} := \{\sum_{k=1}^Kn_k\hat{H}_k(\hat{\theta}_k)^{-1}\}^{-1}\sum_{k=1}^Kn_k(\hat{H}_k(\hat{\theta}_k))^{-1}\hat{\phi}_k$\;
$\hat{\phi}^{WD} := \tilde{\hat{\phi}}^{WD}I(\tilde{\hat{\phi}}^{WD} \in \Phi) + \hat{\phi}^{SaC}I(\tilde{\hat{\phi}}^{WD}\not \in \Phi)$, where  $\hat{\phi}^{SaC}:= N^{-1}\sum_{k=1}^Kn_k\hat{\phi}_k$.
\end{algorithm}

{The Step 4 in the algorithm is necessary since there is no guarantee that after weighting the estimator  $\tilde{\hat{\phi}}^{WD}$ still belongs to the set $\Phi$ as required in Assumption \ref{assumption:compactness}. However the event $\{\tilde{\hat{\phi}}^{WD} \in \Phi\}$ should happen with probability approaching one. Hence,  the $\hat{\phi}^{SAC}I(\tilde{\hat{\phi}}^{WD} \not \in \Phi)$ term is negligible compared with that of  $\hat{\phi}^{WD}I(\tilde{\hat{\phi}}^{WD} \in \Phi)$.
} We need the following assumption in order to establish the MSE bound and asymptotic properties of the proposed WD estimator. 
 
 \begin{assumption} \label{assumption:boundedness}(\textbf{Boundedness})
There exists constants $\rho_{\sigma},c > 0$ such that for $k\geq 1$,
$$ \vertiii{\Sigma_{S,k}(\theta_k^*)}_2\leq \rho_{\sigma}, H_k\succeq cI_{p_1\times p_1}, $$
where $\Sigma_{S,k}(\theta_k) = \mathbb{E}\psi_{\theta_k}(X_{k,1};\theta_k)\psi_{\theta_k}(X_{k,1};\theta_k)^T $.
\end{assumption}
{ By the definition of $H_k(\theta_k)$, we have that
$$\vertiii{H_k}_2 \leq \vertiii{\Psi_{\theta}^{\theta}(\theta_k^*)^{-1}\Sigma_{S,k}(\theta_k^*){\Psi_{\theta}^{\theta}(\theta_k^*)^{-1}}}_2\leq \vertiii{{\Psi_{\theta}^{\theta}(\theta_k^*)^{-1}}}_2^2\vertiii{\Sigma_{S,k}(\theta_k^*)}_2\leq \frac{\rho_{\sigma}}{\rho_{-}^2},$$
which implies $H_k(\theta_k^*)^{-1}\succeq \frac{\rho_{-}^2}{\rho_{\sigma}}I_{p_1\times p_1}$. On the other hand, the above inequality  leads to $\vertiii{{\Psi_{\theta}^{\theta}(\theta_k^*)^{-1}}}_2 \geq \sqrt{\frac{c}{\rho_{\sigma}}}$, and this indicate a finite upper bound for the norm of the Hessian matrix, just as that  assumed in Jordan et al. (2019) \cite{Jordan-2019} and Duan et al. (2021) \cite{Duan-2021}.}

\begin{theorem}\label{thm:risk bound}
Under Assumptions \ref{assumption: relative sample size} -  \ref{assumption: Local strong convexity} and  \ref{assumption:boundedness}, and Assumption \ref{assumption: smoothness} with $v, v_1\geq  2$ , the mean-squared error of the WD estimator $\hat{\phi}^{WD}$ satisfies 
$$\mathbb{E}\|\hat{\phi}^{WD} - \phi^*\|_2^2 \leq \frac{C_1}{nK} + \frac{C_2}{n^2} + \frac{C_3}{n^2K} + \frac{C_4}{n^3} + \frac{C_5K}{n^{\bar{v}}},$$
for 
$n = NK^{-1}$ and $\bar{v}  = min\{v,\frac{v_1}{2}\}$. 
\end{theorem}
{The $v$ and $v_1$ appeared in Assumption 5 {quantify } the moments of the first two orders of the derivatives of the $M$ function and their corresponding Lipschitz functions.}
When the number of data blocks $K = O(n^{min\{1,\frac{\bar{v}-1}{2}\}})$, the convergence rate of MSE of 
$\hat{\phi}^{WD}$ is $\mathcal{O}((nK)^{-1})$, which is the same as the standard full sample estimator. However, when there are too many data blocks such that  $K >> n$, the convergence rate is reduced to $\mathcal{O}(n^{-2})$. 
Furthermore, 
{if the derivatives of the $M$ function and their corresponding Lipschitz functions are heavy-tailed}, say $\bar{v} < 3$,
the convergence rate is further reduced to $\mathcal{O}(Kn^{-\bar{v}})$.
\begin{theorem}\label{thm: asymptotic normality of WD}
Under Assumptions \ref{assumption: relative sample size} -  \ref{assumption: Local strong convexity} and  \ref{assumption:boundedness}, and Assumption \ref{assumption: smoothness} with $v, v_1\geq  2$, if $K = o(n)$,
$$\big(\hat{\phi}^{WD}  - \phi^*)^T\{\sum_{k=1}^Kn_kH_k(\theta_k^*)^{-1}\}\big(\hat{\phi}^{WD}  - \phi^*) \overset{d}{\rightarrow} \chi_{p_1}^2.$$
\end{theorem}

{Although $\{H_k(\theta_k^*)\}_{k=1}^K$ have bounded spectral norms, $\sum_{k=1}^K\frac{n_k}{N}H_k(\theta_k^*)^{-1}$ may not converge to a fixed matrix in presence of heterogeneity. Thus,  we can only obtain the asymptotic normality of the standardized 
$\sqrt{N}\{\sum_{k=1}^K\frac{n_k}{N}H_k(\theta_k^*)^{-1}\}^{1/2}(\hat{\phi}^{DW} - \phi^*)$. This is why Theorem \ref{thm: asymptotic normality of WD} is formulated in a limiting chi-squared  distribution form.} 

{The asymptotic normality implies that we can construct confidence regions for $\phi$  with confidence} level 
$1- \alpha$  as 
\be
\{ \phi \,| \big(\hat{\phi}^{WD}  - \phi)^T\{\sum_{k=1}^Kn_k\hat{H}_k(\hat{\theta}_k)^{-1}\}\big(\hat{\phi}^{WD}  - \phi) \leq \chi^2_{p_1, \alpha}\} \label{eq:WDCI} 
\ee
 after replacing $\sum_{k=1}^Kn_kH_k(\theta_k^*)^{-1}$ with its sample counterpart $\sum_{k=1}^Kn_k\hat{H}_k(\hat{\theta}_k)^{-1}$, where $\chi^2_{p_1, \alpha}$ is the upper $\alpha$ quantile of the $\chi_{p_1}^2$ distribution.}
The block-specific parameter $\lambda_k$ can also  be of interest. Then given the WD estimator of the common parameter $\phi^*$, a question is that whether a more efficient estimator of $\lambda_k^*$ can be obtained. Specifically, we plug in the WD estimator to each data block and re-estimate  $\lambda_k$. The corresponding updated estimator is denoted as $\hat{\lambda}_k^{(2)}$. Actually, the answer is that $\hat{\lambda}_k^{(2)}$ is not necessarily more efficient than $\hat{\lambda}_k$. Due to space limit, more discussions on this aspect are available in Section 1.8 in SM.

\section{Debiased Estimator for diverging K}\label{sec: bias reduction}

It is noted that $K = o(\sqrt{N})$ is required in both Theorems \ref{thm:risk bound} and  \ref{thm: asymptotic normality of WD} to validate the $\mathcal{O}(N^{-1})$ leading order MSE and limiting chi-squared distribution of the WD estimator. The reason is that the bias of the local estimator $\hat{\theta}_k$ is at order $O_p(n_k^{-1})$, which can accumulate across the data blocks by the weighted averaging. This leads to the   bias of $\sqrt{N}(\hat{\phi}^{WD} - \phi^*)$ being at order $O_p(KN^{-1/2})$, which is not necessarily diminishing to zero unless $K = o(\sqrt{N})$. It is worth mentioning that Duan et al. (2021) \cite{Duan-2021} needed the same $K = o(\sqrt{N})$ order in their MLE  framework  to obtain the $\sqrt{N}$-convergence since Li et al. (2003) \cite{Li-2003} showed that the MLE is asymptotically biased when $K/n \rightarrow C \in (0,+\infty)$. This calls for a debias  step for the local estimators before aggregation to  {allow for larger $K$}, {which is needed especially in the Federated Learning scenario where the number of users (data blocks) can be much larger than the size of local data.}

To facilitate the bias correction operation, we have to simplify the notations. 
Suppose $F(\theta)$ is a $p\times 1$ vector function, $\nabla F(\theta)$ is the usual Jacobian whose $l$-th row contains the partial derivatives of the $l$-th element of $F(\theta)$. Then the matrices of higher derivatives are defined recursively so that the $j$-th element of the $l$-th row of $\nabla^sL(\theta)$ (a $p\times p^s$ matrix) is the $1\times p$ vector  $f_{lj}^{v}(\theta) = \partial f_{lj}^{v-1}(\theta)/\partial \theta^T$, where $f_{lj}^{v-1}$ is the $l-$th row and $j$-th element of $\nabla^{v-1}F(\theta)$. We use $\otimes$ to denote a usual Kronecker product. Using Kronecker product we can express $\nabla^vF(\theta) = \frac{\partial^v F(\theta)}{\partial \theta^T \otimes \partial \theta^T \otimes \cdots\otimes\partial \theta^T}$. Besides, let $M_{n,k}(\theta_k) = n_k^{-1}\sum_{i=1}^{n_k}M(X_{k,i};\theta_k)$,
\bea
H_{3,k}(\theta_k) &=& \mathbb{E}\nabla_{\theta_k}^2\psi_{\theta_k}(X_{k,1};\theta_k), \quad 
Q_k(\theta_k) = \{-\mathbb{E}\nabla_{\theta_k}\psi_{\theta_k}(X_{k,1};\theta_k)\}^{-1},\nn\\
\quad d_{i,k}(\theta_k) &=& Q_k(\theta_k)\psi_{\theta_k}(X_{k,i};\theta_k)\quad \hbox{ and } \quad  v_{i,k}(\theta_k) = \nabla_{\theta_k}\psi_{\theta_k}(X_{k,i},\theta_k) - \nabla_{\theta_k}\Psi_{\theta}(\theta_k)   \nn.
\eea
{According to \cite{Rilstone-1996},} the { leading} order bias 
of $\hat{\theta}_k$ is 
\be
Bias(\hat{\theta}_k) = n_k^{-1}Q_k(\theta_k^*)\big(\mathbb{E}v_{i,k}(\theta_k^*)d_{i,k}(\theta_k^*) + \frac{1}{2}H_{3,k}(\theta_k^*)\mathbb{E}\{d_{i,k}(\theta_k^*)\otimes d_{i,k}(\theta_k^*)\}\big). \label{eq: second order bias}
\ee

{Let} $B_k(\theta_k) = Q_k(\theta_k)\big(\mathbb{E}v_{i,k}(\theta_k)d_{i,k}(\theta_k) + \frac{1}{2}H_{3,k}(\theta_k)\mathbb{E}\{d_{i,k}(\theta_k)\otimes d_{i,k}(\theta_k)\}\big)$, whose the first $p_1$ dimension {associated with} $\phi$ are denoted as $B^1_k(\theta_k)$.
The empirical estimator of $B_k(\theta_k)$ is
\be
 \hat{B}_k(\theta_k) = \hat{Q}_k(\theta_k)\big(n_k^{-1}\sum_{i=1}^{n_k}\hat{v}_{i,k}(\theta_k)\hat{d}_{i,k}(\theta_k) + \frac{1}{2}\hat{H}_{3,k}(\theta_k)n_k^{-1}\sum_{i=1}^{n_k}(\hat{d}_{i,k}(\theta_k)\otimes\hat{d}_{i,k}(\theta_k))\big)\label{eq: bias correction formula}
 \ee
 where $\hat{H}_{3,k}(\theta_k) =  n_k^{-1}\sum_{i=1}^{n_k}\nabla_{\theta_k}^2\psi_{\theta_k}(X_{k,i};\theta_k)$, $\hat{Q}_k(\theta_k) = \{- n_k^{-1}\sum_{i=1}^{n_k}\nabla_{\theta_k}\psi_{\theta_k}(X_{k,i};\theta_k)\}^{-1},$ $
 \hat{d}_{i,k}(\theta_k) =\hat{Q}_k(\theta_k)\psi_{\theta_k}(X_{k,i};\theta_k)$ and $\hat{v}_{i,k}(\theta_k) = \nabla_{\theta_k} \psi_{\theta_k}(X_{k,i};\theta_k)$.  
  {Applying bias correction 
to each data block,} 
 we have the bias-corrected local estimator 
\be
\hat{\theta}_{k,bc} := \hat{\theta}_k - n_k^{-1}\hat{B}_k(\hat{\theta}_k)1_{\mathcal{E}_{k,bc}}, \label{eq: debiased local estimator}
\ee
where $\mathcal{E}_{k,bc} = \{\hat{\theta}_k - n_k^{-1}\hat{B}_k(\hat{\theta}_k) \in \Theta_k\}$. The indicator function here is to ensure that $\hat{\theta}_{k,bc}$ is within the parameter space. 

 After the local debiased estimators are obtained, we need to aggregate them with estimated  weights. However, a direct aggregation will invalidate the bias correction procedure due to the correlation between the estimated weights and the {local debiased estimator as they are constructed with the same dataset. The accumulation of the dependence over a large number  of data blocks} can make the bias correction fail.  

To remove such correlation between the local estimators and the corresponding estimated local weights $\hat{W}_k = \{\sum_{s=1}^K \hat{H}_s(\hat{\theta}_s)^{-1}\}^{-1}\hat{H}_k(\hat{\theta}_k)^{-1}$,  
 we first divide each local dataset $\{X_{k,i}\}_{i=1}^{n_k}$ into two equal-sized 
subsets $D_k^{s} = \{X_{k,i}^{(s)}\}_{i=1}^{n_k/2}, s = 1,2$. Then, for $s = 1,2$ we  calculate the local M-estimators $\hat{\theta}_{k,s}$ 
and obtain $\hat{H}_{k,s}(\hat{\theta}_{k,s})$, which is the leading principal  sub-matrix of order $p_1$ of 
\be
(\nabla_{\theta_k}\hat{\Psi}_{\theta_k})^{-1}(\frac{1}{n_k/2}\sum_{i=1}^{n_k/2}\psi_{\theta_k}(X_{k,i}^{(s)};\hat{\theta}_{k,s})\psi_{\theta_k}(X_{k,i}^{(s)};\hat{\theta}_{k,s})^T)(\nabla_{\theta_k}\hat{\Psi}_{\theta_k})^{-T},\nn
\ee
 where $\hat{\Psi}_{\theta_k} = \frac{1}{n_k/2}\sum_{i=1}^{n_k/2}\psi_{\theta_k}(X_{k,i}^{(s)};\hat{\theta}_{k,s})$. We then  perform  the local bias correction
 to  $\{\hat{\theta}_{k,s}\}$ based on data in subset $D_{k}^s$ to attain the debiased estimators $\{\hat{\theta}_{k,s}^{bc}\}$.
 At last, two debiased weighted distributed estimators of the form
 \be
 \tilde{\hat{\phi}}^{dWD}_{s} := \{\sum_{k=1}^Kn_k\hat{H}_{k,s}(\hat{\theta}_{k,s})^{-1}\}^{-1}\sum_{k=1}^Kn_k(\hat{H}_{k,s}(\hat{\theta}_{k,s}))^{-1}\hat{\phi}_{k,2 - |s - 1|}^{bc}\nn
 \ee
 for $s = 1,2$  are averaged to obtain the final debiased WD (dWD) estimator, whose procedure is summarized in Algorithm 2. 
That  the weight estimation and the debiasing  are conducted on different data splits  
 remove the correlation, 
 and realize the gain of bias-correction procedure. 
\begin{algorithm}[!ht]
\caption{\textbf{d}ebiased \textbf{W}eighted \textbf{D}istributed  (dWD) Estimator}
\LinesNumbered
\KwIn{$\{X_{k,i},k = 1,...,K; i = 1,...,n_k\}$}
\KwOut{$\hat{\phi}^{dWD}$}
For each data block, split the data set into two non-overlapping equal-sized subsets and denote those subsets as $D_k^{s} = \{X_{k,i}^{(s)}\}_{i=1}^{n_k/2}, s = 1,2$ \;
Obtain the initial estimates $\hat{\theta}_{k,s} = (\hat{\phi}_{k,s},\hat{\lambda}_{k,s})$ based on data from $D_k^s, s = 1,2$ \;
Calculate $\hat{H}_{k,s}(\hat{\theta}_{k,s})$ in each block $(s = 1,2)$, which is the leading principal  sub-matrix of order $p_1$ of 
$(\nabla_{\theta_k}\hat{\Psi}_{\theta_k}^{(s)})^{-1}(2n_k^{-1}\sum_{i=1}^{n_k/2}\psi_{\theta_k}(X_{k,i}^{(s)};\hat{\theta}_{k,s})\psi_{\theta_k}(X_{k,i}^{(s)};\hat{\theta}_{k,s})^T)(\nabla_{\theta_k}\hat{\Psi}_{\theta_k}^{(s)})^{-T}$ where $\hat{\Psi}_{\theta_k} = 2n_k^{-1}\sum_{i=1}^{n_k/2}\psi_{\theta_k}(X_{k,i}^{(s)};\hat{\theta}_{k,s})$\;
Calculate the bias corrected estimators in each block $(k = 1,2,\cdots,K;s = 1,2)$: $\hat{\theta}_{k,s}^{bc}:=\hat{\theta}_{k,s} - 2n_k^{-1}\hat{B}_{k,s}(\hat{\theta}_{k,s})1_{\mathcal{E}_{k,bc,s}}$ where $\mathcal{E}_{k,bc,s} := \{\hat{\theta}_{k,s} - 2n_k^{-1}\hat{B}_{k,s}(\hat{\theta}_{k,s})\in \Theta_k\}$. Denote the first $p_1$ dimensions of $\hat{\theta}_{k,s}^{bc}$ as $\hat{\phi}_{k,s}^{bc}$\;

Send $\{\hat{\phi}_{k,s}^{bc},\hat{H}_{k,1}(\hat{\theta}_{k,s})^{-1}, s = 1,2\}$ to a central server and construct $\tilde{\hat{\phi}}^{dWD}_{s} := \{\sum_{k=1}^Kn_k\hat{H}_{k,s}(\hat{\theta}_{k,s})^{-1}\}^{-1}\sum_{k=1}^Kn_k(\hat{H}_{k,s}(\hat{\theta}_{k,s}))^{-1}\hat{\phi}_{k,2 - |s - 1|}^{bc}$\;
$\hat{\phi}^{dWD}_s := \tilde{\hat{\phi}}^{dWD}_sI(\tilde{\hat{\phi}}^{dWD}_s \in \Phi) +  K^{-1}\sum_{k=1}^Kn_k\hat{\phi}_{k,2 - |s-1|}^{bc}I(\tilde{\hat{\phi}}^{dWD}_s\not \in \Phi)$ for $s = 1,2$.\;
$\hat{\phi}^{dWD} = \frac{1}{2}\sum_{s=1}^2\hat{\phi}_s^{dWD}.$

\end{algorithm}

To provide theoretical guarantee on the bias correction, we need  an assumption on the third  derivative of the M-function $M$ (see \cite{Zhang-2013}), {which strengthens part of Assumption 5}. 
\begin{assumption}\label{assumption: lip of the third order derivative}(\textbf{Strong smoothness})
For each $x\in \mathbb{R}^p$, the third order derivatives 
of $M(x;\theta_k)$ with respect to $\theta_k$   exist and are $A(x)-$ Lipschitz continuous, i.e.
$$\|(\nabla_{\theta_k}^2\psi_{\theta_k}(x;\theta_k) - \nabla_{\theta_k}^2\psi_{\theta_k}(x;\theta_k^{'}))(u\otimes u)\|_2\leq A(x)\|\theta_k - \theta_k^{'}\|_2\|u\|_2^2,$$
for all $\theta_k,\theta_k^{'}\in U_k$ defined in Assumption \ref{assumption: smoothness} and $u\in \mathbb{R}^p$, where $\mathbb{E}A(X_{k,i})^{2v} \leq A^{2v}$ for some $v>0$ and $A < \infty$.
\end{assumption}

\begin{theorem}\label{thm: MSE bound of dWD}
Under Assumptions \ref{assumption: relative sample size} -  \ref{assumption: Local strong convexity} and \ref{assumption:boundedness} - \ref{assumption: lip of the third order derivative}, and Assumption \ref{assumption: smoothness} with $v, v_1\geq  4$ ,
$$ \mathbb{E}\|\hat{\phi}^{dWD} - \phi^*\|_2^2\leq \frac{C_1}{nK} + \frac{C_2}{n^2K} + \frac{C_3}{n^3} +  \frac{C_4K}{n^{\bar{v}}},$$
where $\bar{v} = min\{v,\frac{v_1}{2}\}$.
\end{theorem} 

{The main difference between the upper bounds} in Theorem \ref{thm: MSE bound of dWD} and that of Theorem \ref{thm:risk bound} for the WD estimator is  the disappearance of the $\mathcal{O}(n^{-2})$ term for the WD estimator, which has been dissolved and absorbed into  the $\mathcal{O}((n^2K)^{-1})$ and $\mathcal{O}(n^{-3})$ terms for the dWD estimator. {As shown  next, this translates to more relaxed  $K=o(n^2)$ as compared with $K=o(n)$ for the WD estimator in Theorem \ref{thm: asymptotic normality of WD}. }

\begin{theorem}\label{thm: Asymptotic normality of dWD}
Under the conditions required by Theorem \ref{thm: MSE bound of dWD}, if $K = o(n^2)$,
$$(\hat{\phi}^{dWD}- \phi^*)^T\{\sum_{k=1}^K n_kH_k(\theta_k^*)^{-1}\}(\hat{\phi}^{dWD}- \phi^*) \overset{d}{\rightarrow}\chi_{p_1}^2.$$
\end{theorem}

Note that the reason why Theorem \ref{thm: Asymptotic normality of dWD} is formulated in the chi-squared distribution form is the same as that when we formulate Theorem \ref{thm: asymptotic normality of WD}, 
 and similar confidence region with confidence level $1-\alpha$ can be constructed as 
\be
\{\phi| \big(\hat{\phi}^{dWD}  - \phi)^T\{\sum_{k=1}^Kn_kH_k(\hat{\theta}_k)^{-1}\}\big(\hat{\phi}^{dWD}  - \phi) \leq \chi_{p_1,\alpha}^2\}.  \label{eq: CI of dWD}
\ee
The fact that the confidence regions of dWD and WD estimators use the same standardizing matrix $\sum_{k=1}^Kn_kH_k(\hat{\theta}_k)^{-1}$ reflects that 
the 
dWD and WD estimators have the same estimation efficiency. However, the debiased version has more relaxed constraint on $K=O(n^2)$ (which is equivalent to $K = o(N^{2/3})$) than that of the WD estimator at $K=o(n)$ ($K = o(\sqrt{N})$).

A more communication-efficient estimator of the common parameter can be defined as the following debiased SaC (dSaC) estimator:
\be
\hat{\phi}^{dSaC} =  N^{-1}\sum_{k=1}^Kn_k(\hat{\phi}_k  - n_k^{-1}\hat{B}_{k}^{1}(\hat{\theta}_k)1_{\mathcal{E}_{k,bc}}),\label{eq: dSaC formula}
\ee
which only performs bias correction and may be preferable when the heterogeneity is not large. The asymptotic property of the dSaC estimator is summarized in the following proposition.
\begin{theorem}\label{thm:  dSaC}Under the conditions required by Theorem \ref{thm: MSE bound of dWD}, if $K = o(n^2)$,
\bea
\mathbb{E}\|\hat{\phi}^{dSaC} - \phi^*\|_2^2 \leq \frac{C_1}{nK} +\frac{C_2}{n^2K} + \frac{C_3}{n^3}\quad \hbox{ and } &&\nn\\
N^2(\hat{\phi}^{dSaC}- \phi^*)^T\{\sum_{k=1}^K n_kH_k(\theta_k^*)\}^{-1}(\hat{\phi}^{dSaC}- \phi^*) \overset{d}{\rightarrow}\chi_{p_1}^2.&&\nn
\eea
\end{theorem}

The corresponding  confidence region with confidence level $1 -\alpha$ can be constructed as \be
\{\phi| N^2\big(\hat{\phi}^{dSaC}  - \phi)^T\{\sum_{k=1}^Kn_kH_k(\hat{\theta}_k)\}^{-1}\big(\hat{\phi}^{dSaC}  - \phi) \leq \chi_{p_1,\alpha}^2\}.\label{eq: CI based on dSaC}
\ee
 It is noted that the dSaC  and SaC estimators have the same asymptotic distribution. Hence, the confidence regions based on the SaC estimator can be constructed as (\ref{eq: CI based on dSaC}) with $\hat{\phi}^{dSaC}$ replaced by $\hat{\phi}^{SaC}$.}  

 To compare with the subsampled average mixture method (SAVGM) estimator proposed in \cite{Zhang-2013} which also performs local bias correction but under the  homogeneous setting, we have the following corollary  to Theorem  \ref{thm:  dSaC}. 
\begin{corollary}\label{corollary: homogeneous}
Under the homogeneous case such that $\{X_{k,i},k = 1,...,K, i = 1,...,n; \}$ are IID distributed, and the assumptions required by Theorem \ref{thm: MSE bound of dWD}, 
\begin{small}
\be
\mathbb{E}\|\hat{\theta}^{dSaC}  - \theta_1^*\|_2^2 \leq \frac{2\mathbb{E}\|\nabla_{\theta_1}\Psi_{\theta}(\theta_1^*)^{-1}\psi_{\theta_1}(X_{1,1};\theta_1^*)\|_2^2}{nK} +\frac{C_1}{n^2K} + \frac{C_2}{n^3},\label{eq: dSaC bound under homogeneity.}
\ee
where $\theta_1^*$ is the true parameter for all the $K$ data blocks.
\end{small}
\end{corollary}
{The SAVGM estimator  resamples $\lfloor rn_k \rfloor $ data points from each data block $k$ for a $r\in(0,1)$ to obtain a local estimator $\hat{\theta}^{SaC}_{k,r}$ based on the sub-samples.}  Then,  the SAVGM estimator is 
\be
\bar{\theta}_{SAVGM} = \frac{\hat{\theta}_k^{SaC} - r\hat{\theta}^{SaC}_{k,r}}{1 - r},\label{eq: SAVGM}
\ee
whose MSE bound as given in Theorem 4 of \cite{Zhang-2013} is 
\be
\mathbb{E}\|\bar{\theta}_{SAVGM}  - \theta_1^*\|_2^2 \leq \frac{2 + 3r}{(1-r)^2}\frac{\mathbb{E}\|\nabla_{\theta_1}\Psi_{\theta}(\theta_1^*)^{-1}\psi_{\theta_1}(X_{1,1};\theta_1^*)\|_2^2}{nK} +\frac{C_1}{n^2K} + \frac{C_2}{n^3}.\label{eq: SAVGM bound}
\ee
Thus, 
the  MSE bound (\ref{eq: SAVGM bound}) of the SAVGM estimator  has an inflated factor $\frac{2 + 3r}{2(1 -r)^2} > 1$ for $r \in (0,1)$, when compared with that of the dSaC estimator, { although it   
 is computationally more efficient than the  dSaC and dWD estimators as it only draws one subsample in its resampling. } For more comparisons between the dSaC estimator and  one-step estimators proposed by Huang and Huo (2019) \cite{Huang-Huo-2019}, see Section 1.10 in SM.
\section{Simulation Results}\label{sec: Simulation}

We report results from simulation experiments 
 designed to verify two sets  theoretical findings made in the previous sections.   One was to confirm 
 the finding in Section \ref{sec: Preliminaries} that the full sample estimator $\hat{\phi}_{full}$ is not necessarily more efficient than the {SaC} estimator $\hat{\phi}^{SaC}$. The other was to evaluate  the numerical performance of the newly proposed { weighted distributive ({WD}) , debaised SaC ({dSaC}) and debiased WD ({dWD})}  estimators of the common parameter 
 and compare them with the existing {SaC} and {subsampled average mixture method (SAVGM)} 
 (with subsampling rate $r = 0.05$) estimators. 
{
Although the SAVGM estimator \citep{Zhang-2013} 
was proposed under the homogeneous setting, but since its main bias correction is performed locally on each data block $k$ as shown in  (\ref{eq: SAVGM}), similar theoretical bounds as formula (\ref{eq: SAVGM bound}) can be derived without much  modifications on the original proof.}  
  Throughout the simulation experiments, the results of each simulation setting were based on $B = 500$ number of replications and were conducted in R paralleled with a single 10-core Intel(R) Core(TM) i9-10900K @3.7 GHz processor.

In the first simulation experiment, we simulated  the {errors-in-variables} Model (\ref{eq:error-in-variable}) with the objective function (\ref{eq: obj for EIV model}) to compare the  performance of { the full sample, the SaC and the WD estimators}:  $\hat{\phi}_{full}$,  $\hat{\phi}^{SaC}$ and  $\hat{\phi}^{WD}$. The simulation was carried out by first generating IID $\{Z_{i,k}\}$ from $\mathcal{N}(\mu_Z,\sigma^2_Z)$, and then 
upon given a $Z_{i,k}$,  
$(X_{k,i},Y_{i,k})^T$ were independently drawn from $\mathcal{N}\big((Z_{i,k},\phi^* + \lambda_k^*Z_{i,k})^T,\sigma^2I_{2\times2}\big)$. 
We chose $\phi^* = 1,K = 2,\sigma^2 = 1$ and $n_1 = n_2 = 5 \times 10^4 = N / 2$, and  $\lambda_1^*, \lambda_2^*, \mu_Z$ and $ \sigma^2_Z$  were those reported  in Table \ref{tb: estimated bias and sd} under four scenarios.  

As discussed in Section \ref{sec: Preliminaries}, the  relative efficiency of $\hat{\phi}_{full}$ to  $\hat{\phi}^{SaC}$ depends on the ratio \\ $\sigma^2(\mathbb{E}Z)^2 / (var(Z)\mathbb{E}Z^2)$ as shown in 
(\ref{eq: compare SaC and full in EIV}). 
We designed four scenarios according to the above ratio under $\lambda_1^*\not = \lambda_2^*$ and $\mathbb{E}Z \not = 0$, {respectively}, which represented the settings where the full sample estimator $\hat{\phi}_{full}$ would be less  (Scenario 1) or more (Scenario 2) efficient than the SaC estimator as predicted by the ratio, but not as efficient as the weighted distributed estimator $\hat{\phi}^{WD}$. Scenario 3 ($\lambda_1^*\not = \lambda_2^*, \mathbb{E}Z = 0$) was the case when $\hat{\phi}_{full}$ and $\hat{\phi}^{WD}$ would be asymptotically equivalent, and both estimators would be more efficient than $\hat{\phi}^{SaC}$.   Scenario  4 was the homogeneous case with $\lambda_1^* = \lambda_2^*$ in which all the three estimators would have the same asymptotic efficiency.  
%
{For all the four scenarios, the ARE column of the Table \ref{tb: estimated bias and sd} confirmed the relative efficiency as predicted by 
the asymptotic variances in (\ref{eq: compare SaC and full in EIV}), and  was well reflected in the comparison of the RMSEs}, as the bias is of smaller order as compared with that of the SD and thus negligible.

\begin{table}[!ht]
\centering

\caption{Average root mean squared error (RMSE) and the standard
deviation (SD), multiplied by $10^2$,  of the full sample estimator $\hat{\phi}_{full}$, the SaC estimator $\hat{\phi}^{SaC}$ and the WD estimator $\hat{\phi}^{WD}$ under  four scenarios  for the errors-in-variables model (12)  for $N=10^5, K=2$ and $n_1 =n_2$. AREs  (asymptotic relative efficiency)     
 of $\hat{\phi}_{full}$ to $\hat{\phi}^{SaC}$ are calculated from  (\ref{eq: compare SaC and full in EIV}). }
\begin{tabular*}{\textwidth}{ccc@{\extracolsep{\fill}}cccccc}
\toprule
&  & &\multicolumn{2}{c}{$\hat{\phi}_{full}$} & \multicolumn{2}{c}{$\hat{\phi}^{SaC}$}  & \multicolumn{2}{c}{$\hat{\phi}^{WD}$}\\
Scenario & $(\lambda_1^*,\lambda_2^*)$  & ARE & RMSE & SD & RMSE & SD  & RMSE & SD\\ 
\midrule
Scenario 1  & (0.25,3.25) & 0.89 & 4.55  & 4.51  & 4.12  & 4.09 & 3.91 & 3.89 \\
$(\mu_Z = 1,\sigma_Z^2 = 0.1)$ & (0.5,3.5) & 0.93 & 4.65& 4.65 & 4.35 & 4.35 & 4.08 & 4.08\\
& (0.75,3.75) & 0.97 & 4.52& 4.52 &  4.40& 4.38& 4.13 & 4.13\\
\hline
Scenario 2  & (0.25,2.25) & 1.18 & 2.95& 2.95 & 3.24 & 3.24 & 2.89 & 2.89 \\
$(\mu_Z = 3,\sigma_Z^2 = 0.5)$ & (0.75,2.75) & 1.28 & 3.28 & 3.26 & 3.65 & 3.64 & 3.17 & 3.16\\
& (1.25,3.25) & 1.31& 3.71 & 3.71 & 4.16 & 4.07 & 3.64 &  3.61 \\
\hline
Scenario 3 & (0.25,2.25) & 1.97 & 0.41 & 0.41 & 0.61 & 0.61 & 0.41 & 0.41 \\
$(\mu_Z = 0,\sigma_Z^2 = 0.5)$ & (0.75,2.75) & 1.92 &0.51 & 0.51 & 0.70 & 0.70 & 0.51 & 0.51 \\
& (1.25,3.25) & 1.68 &0.64 & 0.64 & 0.82 & 0.82 & 0.64 & 0.64\\
\hline
Scenario 4  & (0.5,0.5) &1 & 3.25 & 3.24&  3.31& 3.28 & 3.30 & 3.26 \\
$(\mu_Z = 4,\sigma_Z^2 = 0.5)$ & (1.0,1.0) &1 & 3.53& 3.53& 3.59 & 3.59 & 3.59 & 3.59\\
& (1.5,1.5) & 1& 4.06& 4.03 & 4.08 & 4.07 & 4.06 & 4.06\\
\bottomrule
\end{tabular*}
\label{tb: estimated bias and sd}
\end{table}

In the second simulation experiment, we evaluated the numerical performance of the five estimators 
{for the common parameter $\phi$} under a  logistic regression model. 
 For each of $K$ data block with $K \in \{10,50,100,250, 500,1000,2000\}$,  
$\{(X_{k,i};Y_{k,i})\}_{i=1}^n\subset \mathbb{R}^{p}\times\{0,1\}$ were independently sampled from the following model: 
$$X_{k,i}\overset{i.i.d}{\sim}\mathcal{N}(\mathbf{0}_{p\times 1},0.75^2I_{p\times p}) \quad \hbox{and} \quad P(Y_{k,i} = 1 |X_{k,i}) = \frac{exp( X_{k,i}^T\theta_k^*)}{1 + exp( X_{k,i}^T\theta_k^*)},$$
where $\theta_k^* = (\phi^*,\lambda_k^{*T})^T$, $\phi^* = 1$, $\lambda_{k}^* = (\lambda_{k,1}^*,\lambda_{k,2}^*,\cdots,\lambda_{k,p_2}^*)^T$ and $\lambda_{k,j}^* = (-1)^j 10(1 - \frac{2(k-1)}{K-1}) $. 
 The sample sizes of the data blocks were equal at $n  = NK^{-1}$ with $N = 2\times 10^6$.   Two levels of  the dimension $p_2=4$ and $10$ of the nuisance parameter $\lambda_k$ were considered.  
A derivation of the bias correction formula {for the logistic model}  is given in Section 1.9 of the SM. 

\begin{figure}[!ht]
    \centering
    \begin{supertabular}{cc}
    
         (a) Absolute Bias ($p_2 = 4$) &  (b) RMSE ($p_2 = 4$) 
         \\
         \includegraphics[width = 7.2cm,height = 5.4cm]{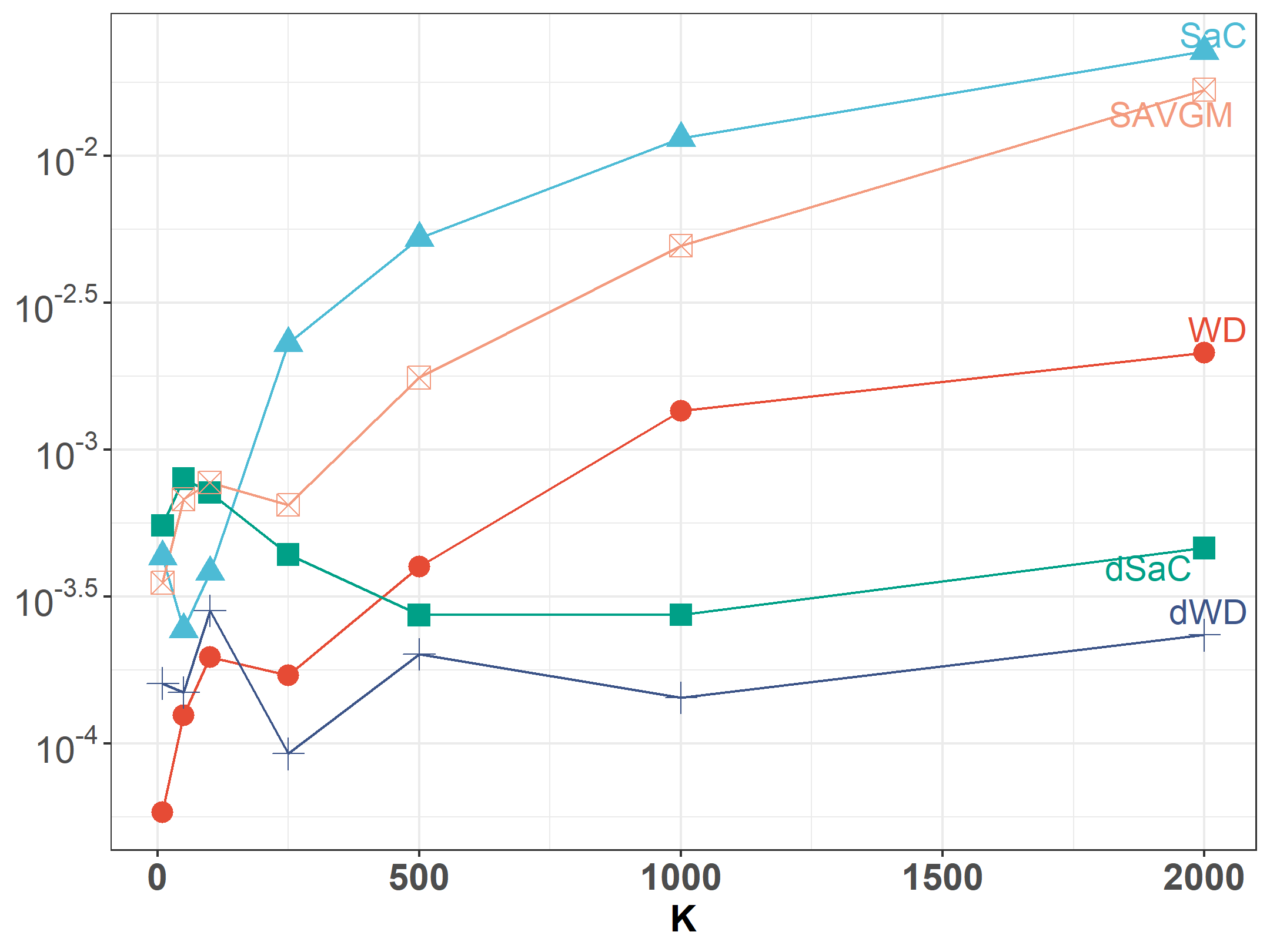}&
          \includegraphics[width = 7.2cm,height = 5.4cm]{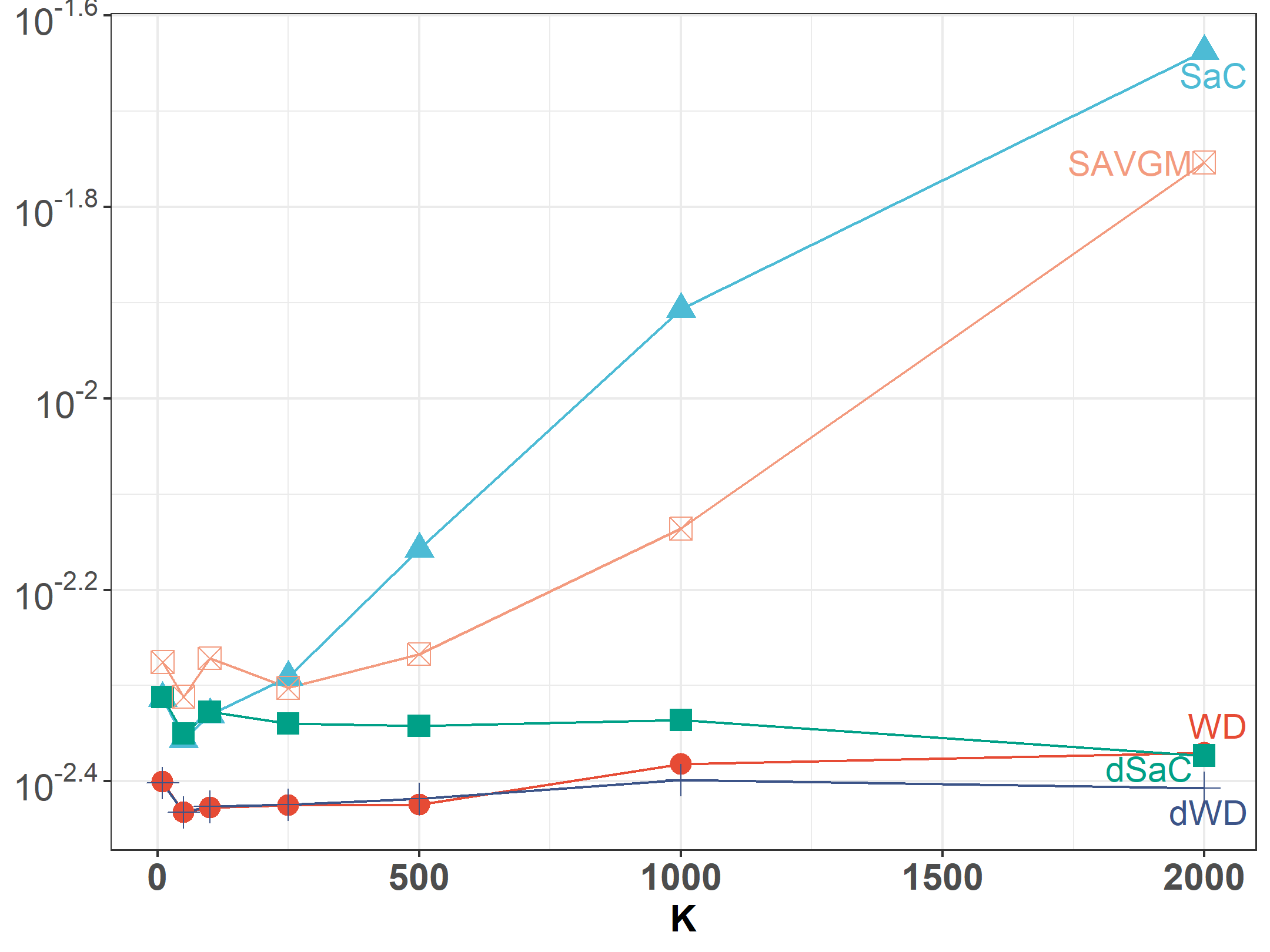}
          \\
          (c) Absolute Bias ($p_2 = 10$)&  (d) RMSE ($p_2 = 10$) 
         \\
         \includegraphics[width = 7.2cm,height = 5.4cm]{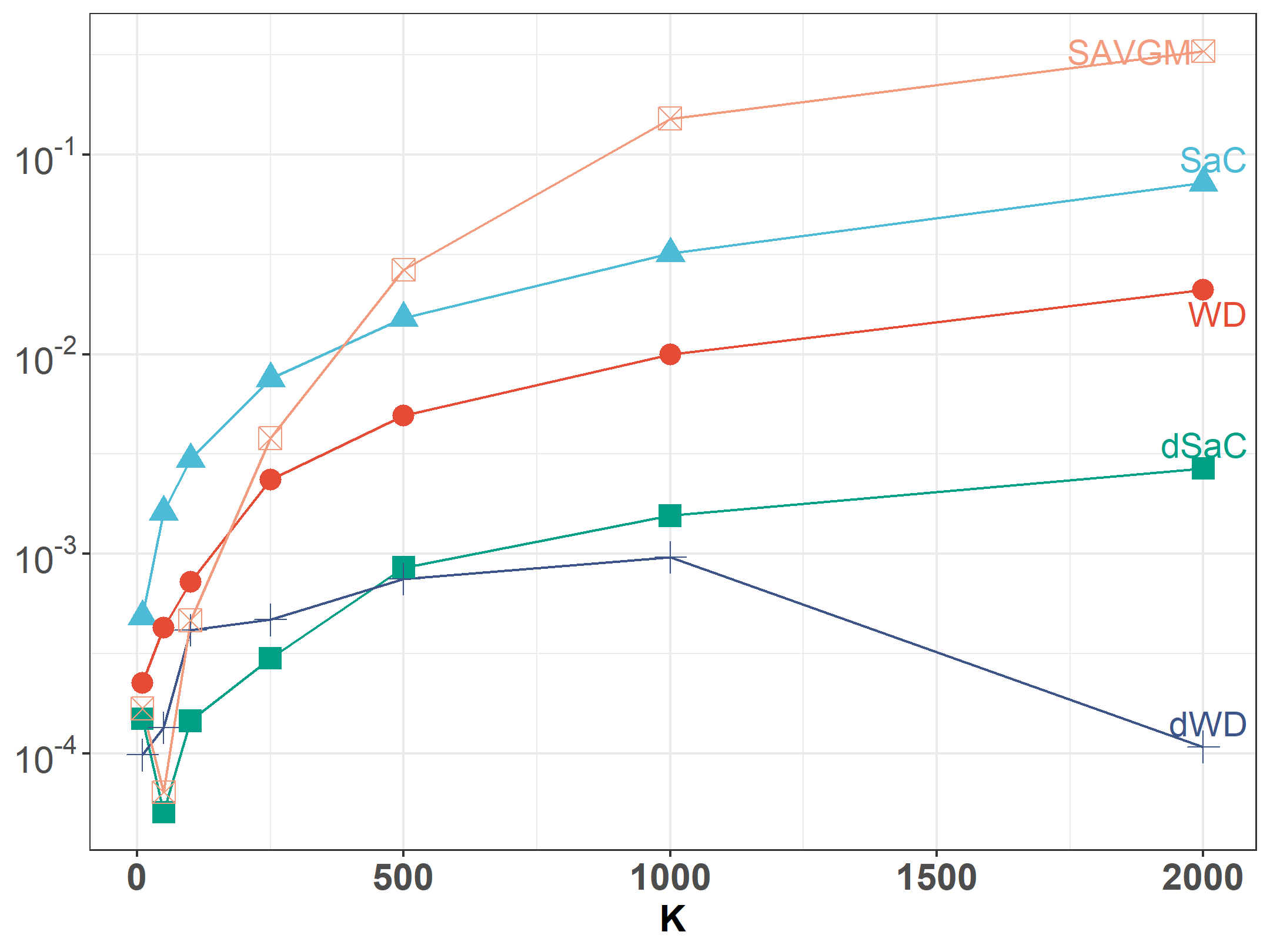}&
          \includegraphics[width = 7.2cm,height = 5.4cm]{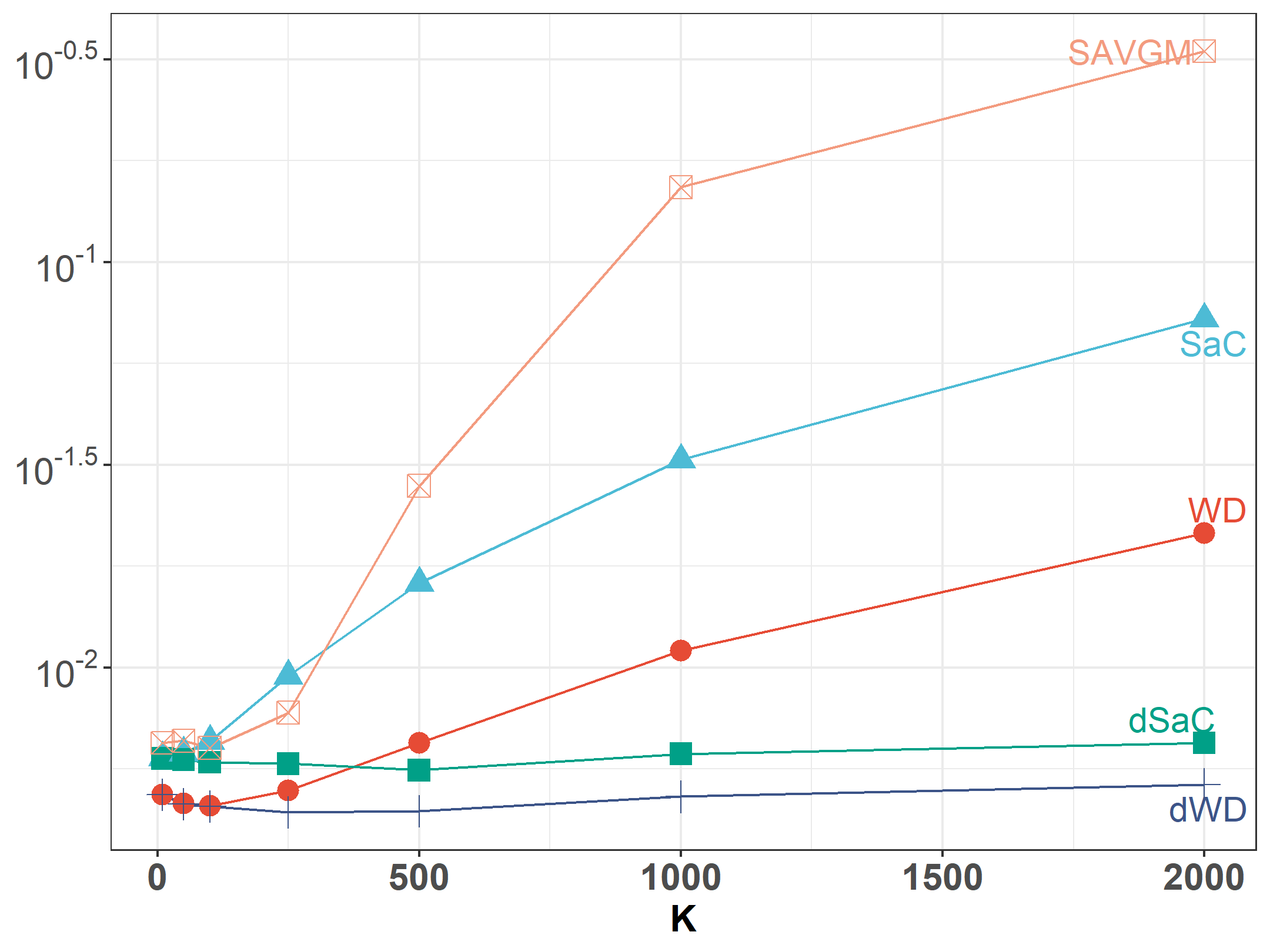}
          \\
    \end{supertabular}
    \caption{Average simulated 
    bias (a, c) and the root mean square errors (RMSE) (b,d)  of the weighted distributed  (WD) (red circle), 
    the SaC (blue triangle), 
    the debiased SaC (dSaC) (green square), 
    the debiased WD (dWD) (purple cross), 
    the subsampled average mixture SAVGM (pink square cross)  estimators,  with respect to the number of data block  $K$  for the logistic regression model with the dimension $p_2$ of the  nuisance parameter $\lambda_k$ being $4$ and $10$, respectively with the full  sample size $N = 2\times 10^6$. 
    }
    \label{fig:fix K}
\end{figure}

 Figure \ref{fig:fix K} reports {the  root mean square errors (RMSEs)
 and absolute bias of the estimators}. 
It is observed  that the weighted distributed estimator WD, and the two debiased estimators dSaC and dWD had 
 smaller RMSE than those of   the 
 SaC and SAVGM for almost all the simulation settings.   Between the SaC  and SAVGM, the SAVGM fared better in the lower dimensional case of $p_2=4$, but was another way around for $p_2=10$.  
It was evidence that the WD estimator had much smaller RMSEs than the SaC and SAVGM estimators for all the block number $K$, realizing its theoretical promises.   { In most cases the WD estimator had smaller bias than the SaC estimator although it was not debiased.} It also had smaller RMSEs than the debiased SaC estimator dSaC for almost all cases of the block numbers for $p_2=4$, while in the higher dimensional $p_2=10$ the WD estimator was advantageous for $K \le 250$. The latter indicated  the need for conducting the bias correction to the WD estimator. 
Both bias corrected  dWD and dSaC were very effective in reducing the bias of the WD and SaC estimators, respectively, especially for larger $K$ when the bias was more severe.  The debiased WD 
attained the smallest RMSEs and the bias in all settings, suggesting the need for conducting  both weighting and the bias correction in the distributed inference especially for large $K$. 
These empirical results were consistent with  Theorems \ref{thm:risk bound} and \ref{thm: MSE bound of dWD}, namely  the leading RMSE term of  the WD estimator changes from $\mathcal{O}((Kn)^{-1})$ to $\mathcal{O}(n^{-2})$ when 
  $K$ surpasses the local sample size $n$, while the leading RMSEs of the dWD is still $\mathcal{O}((nK)^{-1})$ until $K >> n^2$.

We also evaluated the coverage probabilities and widths of the $1-\alpha$ ($\alpha = 0.01,0.05,0.1$) confidence intervals (CIs) of the common parameter based on the asymptotic normality 
as  
 given after Theorems \ref{thm: asymptotic normality of WD} and \ref{thm: Asymptotic normality of dWD}. 
The SAVGM estimator was not included as its asymptotic distribution was not made available in \cite{Zhang-2013}. 
Table \ref{tb:CIs for different p2} reports the empirical coverage and the average width of the CIs.  
{It is observed that for the lower dimensional nuisance parameter case of  $p_2=4$  the four types of the  CIs  all had quite adequate coverage levels when $K \le 100$.  However, for $K \ge 250$,  the SaC CIs first started to lose coverage, followed by those of the  WD, while the CIs of the debiased SaC (dSaC) and debiased weighted distributed  (dWD) estimators can hold up to the promised coverage for all cases of $K$.  The outstanding performance of the dSaC and dWD CIs was largely replicated for the higher dimensional nuisance parameter case of $p_2=10$, while the other two non-debiased estimator based CIs had their coverage quickly slipped below the nominal coverage levels. Although the dSaC CIs had comparable coverages with the dWD CIs, their widths were much wider than those of the dWD. This was largely due to the fact that the weighted averaging conducted in the weighted distributed estimation  reduced the variation and hence the width of the CIs. 
The widths of the WD CIs were largely the same with those of the dWD, and yet the coverage levels of the dWD CIs were much more accurate indicating the importance of the bias correction as it shifted the CIs without inflating the width. } 

In addition to the simulation experiments on the statistical properties of the estimators, the computation efficiency of the estimators was also evaluated. 
Table \ref{tb: all time} 
reports the {average} CPU time per simulation run based on 500 replications of the five estimators for a range of  $K$ and dimension  $p_2$ of the nuisance parameter for the logistic regression model with 
the total sample size $ N = 2\times 10^6$.

The computation speed of the dSaC and dWD estimators were relatively slower than those of  the SaC, WD  and SAVGM estimators. 
The WD estimator was quite fast, which means that the re-weighting { used less computing time}  than the bias-reduction. In comparison, the dWD estimator was the slowest as a cost for attaining the best RMSE among the five estimators in all settings. 
It is observed in Table \ref{tb: all time}  that the overall computation time for   each estimator first decreased and then increased as $K$ became larger. The decrease in time was because the benefit of the distributed computation, while the increase { was due to the increase in the number of optimization associated with the M-estimation performed as $K$ got larger.}  However, it is worth mentioning that these results did not account for  the potential time expenditure in data communication among different data blocks.

\section{Discussion} 
\label{sec: Discussion}
This paper investigates several distributed M-estimators in the presence of heterogeneous distributions  among the data blocks. 
The weighted distributed (WD) estimator is able to improve the estimation efficiency of the  "Split-And-Conquer" (SaC) estimator for the common parameter. 
Two debiased estimators ( dWD and dSaC) 
are proposed to allow for larger numbers of data blocks $K$. 
The statistical properties of these three estimators are shown to be advantageous over the SaC and  SAVGM estimators. 
In particular, the WD estimator has good performance  for smaller $K$ relative to $n$, and the debiased WD estimator that conducted both bias correction and weighting offers good estimation accuracy for large $K$. 

An important issue for the distributed estimation is the size of $K$ relative to the local average sample size $n$. This is especially true in Federated Learning setting where the number of clients (data blocks) are usually very large. { Both  SaC and WD estimators}  require $K = o(\sqrt{N})$ 
to preserve the $\mathcal{O}(N^{-1})$ convergence rate for its MSE and the $\sqrt{N}$ rate for the asymptotic variance. {The debiased  dWD and dSaC} relax  the restriction to $K = o(N^{2/3})$ without compromising the convergence rate. 
 The dSaC may be used as a computationally cheaper version of the dWD at the cost of larger variations and wider confidence regions  when compared with dWD.


\bibliographystyle{apalike}
\bibliography{ref}

\setlength\tabcolsep{2.0pt}
\begin{table}[!ht]
\centering
\caption{Coverage probabilities and widths (in parentheses, multiplied by $100$) of the $1-\alpha$  confidence intervals { for the common parameter $\phi$ in the logistic regression model based on the asymptotic normality of} the SaC,  
the WD, 
the debiased SaC 
and the debiased WD  
estimators with respect to the number of data blocks $K$  { for two dimensions of the heterogeneous parameter $p_2$} with the full sample size $ N = 2\times 10^6$.}

\centerline{(a) $p_2 =4$} 
\begin{tabular*}{\textwidth}{cc@{\extracolsep{\fill}}ccc|ccc|ccc|ccc}
\hline
K && \multicolumn{3}{c|}{SaC} &\multicolumn{3}{c|}{WD} & \multicolumn{3}{c|}{dSaC} & \multicolumn{3}{c}{dWD}\\
&$1-\alpha$ & 0.99 & 0.95 & 0.90 &0.99 & 0.95 & 0.90 &0.99 & 0.95 & 0.90 &0.99 & 0.95 & 0.90 \\
\hline
10&&0.99 & 0.96 & 0.92 & 0.99 & 0.97 & 0.91 & 0.99 & 0.96 & 0.92 & 0.99 & 0.96 & 0.91 \\ 
  &&(2.45) & (1.87) & (1.57) & (2.03) & (1.55) & (1.30) & (2.45) & (1.87) & (1.57) & (2.03) & (1.55) & (1.30) \\ 
50&&0.99 & 0.95 & 0.91 & 0.98 & 0.93 & 0.89 & 0.99 & 0.95 & 0.91 & 0.99 & 0.93 & 0.88 \\ 
  &&(2.36) & (1.80) & (1.51) & (1.97) & (1.50) & (1.26) & (2.36) & (1.80) & (1.51) & (1.97) & (1.50) & (1.26) \\ 
100&&0.98 & 0.94 & 0.91 & 0.99 & 0.95 & 0.91 & 0.99 & 0.95 & 0.91 & 0.99 & 0.95 & 0.91 \\ 
  &&(2.36) & (1.79) & (1.51) & (1.96) & (1.49) & (1.25) & (2.36) & (1.79) & (1.51) & (1.96) & (1.49) & (1.25) \\ 
250&&0.99 & 0.93 & 0.85 & 0.99 & 0.95 & 0.90 & 0.99 & 0.96 & 0.91 & 0.99 & 0.95 & 0.90 \\    &&(2.36) & (1.79) & (1.50) & (1.96) & (1.49) & (1.25) & (2.36) & (1.79) & (1.50) & (1.96) & (1.49) & (1.25) \\ 
  500&&0.91 & 0.77 & 0.66 & 0.99 & 0.95 & 0.88 & 0.99 & 0.96 & 0.90 & 0.99 & 0.95 & 0.89 \\ 
  &&(2.36) & (1.80) & (1.51) & (1.96) & (1.49) & (1.25) & (2.36) & (1.80) & (1.51) & (1.96) & (1.49) & (1.25) \\ 
  1000&&0.65 & 0.41 & 0.28 & 0.99 & 0.94 & 0.88 & 0.99 & 0.94 & 0.88 & 0.99 & 0.93 & 0.88 \\ 
  &&(2.38) & (1.81) & (1.52) & (1.96) & (1.49) & (1.25) & (2.38) & (1.81) & (1.52) & (1.97) & (1.50)& (1.25) \\ 
  2000&&0.01 & 0.01 & 0.00 & 0.99 & 0.91 & 0.81 & 0.98 & 0.94 & 0.88 & 0.99 & 0.94 & 0.90 \\ 
  &&(2.42) & (1.84) & (1.55) & (1.96) & (1.50) & (1.25) & (2.42) & (1.84) & (1.55) & (1.98) & (1.50) & (1.26) \\ 
\hline 
\end{tabular*}
\label{tb:CIs for different p2}
\end{table}
\begin{table}[!ht]
\centering 
\centerline{(b) $p_2 = 10$ } 
\begin{tabular*}{\textwidth}{cc@{\extracolsep{\fill}}ccc|ccc|ccc|ccc}
\hline
K && \multicolumn{3}{c|}{SaC} &\multicolumn{3}{c|}{WD} & \multicolumn{3}{c|}{dSaC} & \multicolumn{3}{c}{dWD}\\
&$1-\alpha$ & 0.99 & 0.95 & 0.90 &0.99 & 0.95 & 0.90 &0.99 & 0.95 & 0.90 &0.99 & 0.95 & 0.90 \\
\hline
10 & &0.99 & 0.94 & 0.88 & 1.00 & 0.96 & 0.92 & 1.00 & 0.94 & 0.88 & 1.00 & 0.96 & 0.92 \\ 
  &&(3.05) & (2.32) & (1.95) & (2.41) & (1.84) & (1.54) & (3.05) & (2.32) & (1.95) & (2.42) & (1.84) & (1.54) \\ 
 50&& 0.99 & 0.93 & 0.87 & 0.99 & 0.95 & 0.88 & 0.98 & 0.94 & 0.88 & 0.99 & 0.96 & 0.88 \\ 
  &&(2.94) & (2.24) & (1.88) & (2.29) & (1.74) & (1.46) & (2.94) & (2.24) & (1.88) & (2.29) & (1.74) & (1.46) \\ 
  100&&0.97 & 0.89 & 0.84 & 0.97 & 0.93 & 0.87 & 0.98 & 0.95 & 0.90 & 0.98 & 0.94 & 0.89 \\ 
  &&(2.93) & (2.23) & (1.87) & (2.28) & (1.74) & (1.46)& (2.93) & (2.23) & (1.87) & (2.29) & (1.74) & (1.46) \\ 
  250&&0.89 & 0.72 & 0.63 & 0.98 & 0.92 & 0.87 & 1.00 & 0.97 & 0.90 & 1.00 & 0.96 & 0.90 \\ 
  &&(2.94) & (2.24) & (1.88) & (2.28) & (1.74) & (1.46) & (2.94) & (2.24) & (1.88) & (2.29) & (1.74) & (1.46) \\ 
  500&&0.51 & 0.28 & 0.18 & 0.93 & 0.81 & 0.70 & 0.99 & 0.94 & 0.90 & 0.98 & 0.94 & 0.88 \\ 
  &&(2.97) & (2.26) & (1.90) & (2.29) & (1.74) & (1.46) & (2.97)& (2.26) & (1.90) & (2.30) & (1.75) & (1.47) \\ 
  1000&&0.00 & 0.00 & 0.00 & 0.66 & 0.37 & 0.28 & 0.99 & 0.95 & 0.90 & 0.99 & 0.96 & 0.89 \\ 
  &&(3.04) & (2.31) & (1.94) & (2.30) & (1.75) & (1.47) & (3.04) & (2.31) & (1.94) & (2.34) & (1.78) & (1.49) \\ 
  2000&&0.00 & 0.00 & 0.00 & 0.02 & 0.00 & 0.00 & 0.99 & 0.96 & 0.90 & 0.99 & 0.93 & 0.87 \\ 
  &&(3.22) & (2.45) & (2.06) & (2.34) & (1.78) & (1.49) &(3.22) & (2.45) & (2.06) & (2.40) & (1.82) & (1.53) \\
\hline
\end{tabular*}
\label{tb: CIs for p2 = 10}
\end{table}

\begin{table}[!ht]
\renewcommand\arraystretch{0.85}
\centering
\caption{Average CPU time   for each replication based on  $B = 500$ replications for the SaC,   
 the SAVGM,  
 the WD,   
 the debiased SaC  
 and the debiased WD estimators  
 for the logistic regression model with respect to $K$ and the dimension $p_2$ of the nuisance parameter. Total sample size $ N = 2\times 10^6$. }
\begin{tabular*}{\textwidth}{c@{\extracolsep{\fill}}ccccc}
\toprule
K & $SaC$ & $SAVGM$ & $WD$ & $dSaC$ & $dWD$ \\
\midrule
&& \multicolumn{3}{c}{$p_2 = 4$}\\
\cmidrule(r){2-6}
10 & 15.65 &  15.97 & 18.50 & 20.00 & 21.95 \\ 
   50 & 9.63 & 9.95 & 10.66 & 12.37 & 14.59 \\ 
  100 & 8.09 & 8.63 & 8.76 & 10.50 & 12.05 \\ 
  250 & 8.49 & 9.69 & 9.07 & 10.84 & 12.82 \\ 
  500 & 9.68 & 11.58 & 10.25 & 11.97 & 14.84 \\ 
  1000 & 11.67 & 13.81 & 12.32 & 13.93 & 19.08 \\ 
  2000 & 15.78 & 19.68 & 16.57 & 18.11 & 28.55 \\ 
  && \multicolumn{3}{c}{$p_2 = 10$}\\
\cmidrule(r){2-6}
   10 & 34.60 & 35.19 & 43.84 & 50.47 & 55.35 \\ 
   50 & 20.13 & 20.18 & 24.16 & 29.99 & 33.69 \\ 
  100 & 15.60 & 16.20 & 17.74 & 23.63 & 24.47 \\ 
  250 & 10.77 & 12.61 & 11.88 & 18.22 & 20.39 \\ 
  500 & 11.55 & 14.50 & 12.56 & 18.80 & 23.73 \\ 
  1000 & 15.23 & 18.27 & 16.28 & 22.38 & 32.24 \\ 
  2000 & 23.42 & 27.99 & 24.62 & 30.43 & 48.05 \\ 

\bottomrule
\end{tabular*}
\label{tb: all time}
\end{table}


\end{document}